# MULTISOURCE BAYESIAN SEQUENTIAL CHANGE DETECTION

By Savas Dayanik,[1] H. Vincent Poor[2] and Semih O. Sezer[2]

*Princeton University, Princeton University and University of Michigan*

Suppose that local characteristics of several independent compound Poisson and Wiener processes change suddenly and simultaneously at some unobservable disorder time. The problem is to detect the disorder time as quickly as possible after it happens and minimize the rate of false alarms at the same time. These problems arise, for example, from managing product quality in manufacturing systems and preventing the spread of infectious diseases. The promptness and accuracy of detection rules improve greatly if multiple independent information sources are available. Earlier work on sequential change detection in continuous time does not provide optimal rules for situations in which several marked count data and continuously changing signals are simultaneously observable. In this paper, optimal Bayesian sequential detection rules are developed for such problems when the marked count data is in the form of independent compound Poisson processes, and the continuously changing signals form a multidimensional Wiener process. An auxiliary optimal stopping problem for a jump-diffusion process is solved by transforming it first into a sequence of optimal stopping problems for a pure diffusion by means of a jump operator. This method is new and can be very useful in other applications as well, because it allows the use of the powerful optimal stopping theory for diffusions.

**1. Introduction.** Suppose that at some unobservable disorder time $\Theta$, the local characteristics of several independent compound Poisson and Wiener processes undergo a sudden and simultaneous change. More precisely, the

Received November 2006; revised July 2007.

[1]Supported in part by the Air Force Office of Scientific Research, Grant AFOSR-FA9550-06-1-0496, and by the U.S. Department of Homeland Security through the Center for Dynamic Data Analysis for Homeland Security administered through ONR Grant N00014-07-1-0150 to Rutgers University.

[2]Supported by the U.S. Army Pantheon Program.

pairs $(\lambda_0^{(i)}, \nu_0^{(i)})$, $1 \le i \le m$, consisting of the arrival rate and mark distribution of $m$ compound Poisson processes $(T_n^{(i)}, Z_n^{(i)})_{n\ge 1}$, $1 \le i \le m$, become $(\lambda_1^{(i)}, \nu_1^{(i)})$, $1 \le i \le m$, and $d$ Wiener processes $W_t^{(j)}$, $1 \le j \le d$ gain drifts $\mu^{(j)}$, $1 \le j \le d$ at time $\Theta$.

We assume that $\Theta$ is a random variable with the zero-modified exponential distribution

$$\mathbb{P}\{\Theta = 0\} = \pi \quad \text{and} \quad \mathbb{P}\{\Theta > t\} = (1-\pi)e^{-\lambda t}, \qquad t \ge 0, \tag{1.1}$$

and $(\lambda_0^{(i)}, \nu_0^{(i)})_{1\le i\le m}$, $(\lambda_1^{(i)}, \nu_1^{(i)})_{1\le i\le m}$, $(\mu^{(j)})_{1\le j\le d}$, $\pi$, and $\lambda$ are known. The objective is to detect the disorder time $\Theta$ as soon as possible after disorder happens by using the observations of $(T_n^{(i)}, Z_n^{(i)})_{n\ge 1}$, $1 \le i \le m$, and

$$X_t^{(j)} = X_0^{(j)} + \mu^{(j)}(t-\Theta)^+ + W_t^{(j)}, \qquad t \ge 0, 1 \le j \le d.$$

More precisely, if $\mathbb{F} = \{\mathcal{F}_t\}_{t\ge 0}$ denotes the observation filtration, then we would like to find, if it exists, an $\mathbb{F}$-stopping time $\tau$ whose Bayes risk

$$R_\tau(\pi) \triangleq \mathbb{P}\{\tau < \Theta\} + c\mathbb{E}(\tau - \Theta)^+, \qquad 0 \le \pi < 1 \tag{1.2}$$

is the smallest for any given constant cost parameter $c > 0$ and calculate its Bayes risk. If such a stopping time exists, then it provides the best trade-off between false alarm frequency $\mathbb{P}\{\tau < \Theta\}$ and expected detection delay cost $c\mathbb{E}(\tau - \Theta)^+$.

Important applications of this problem are the quickest detection of manufacturing defects during product quality assurance, online fault detection and identification for condition-based equipment maintenance, prompt detection of shifts in the riskiness of various financial instruments, early detection of the onset of an epidemic to protect public health, quickest detection of a threat to homeland security, and online detection of unauthorized access to privileged resources in the fight against fraud. In many of those applications, a range of data, changing over time either continuously or by jumps or both, are collected from multiple sources/sensors in order to detect a sudden unobserved change as quickly as possible after it happens, and the problems can be modeled as the quickest detection of a change in the local characteristics of several Wiener and compound Poisson processes. For example, in condition-based maintenance, an equipment is monitored continuously by a web of sensors for both continuously-changing data (such as oil level, temperature, pressure) and marked count data (e.g., number, size and type of wear particles in the oil); see Byington and Garga [6]. For the assessment of financial risks of an electricity delivery contract, the spot price of electricity is sometimes modeled by a jump-diffusion process; see, for example, Weron, Bierbrauer and Trück [18] and Cartea and Figueroa [7].

In the past, the Bayesian sequential change-detection problems have been studied for Wiener processes by Shiryaev [17, Chapter 4] and for Poisson processes by Peskir and Shiryaev [14, 15], Gapeev [10], Bayraktar, Dayanik and Karatzas [2, 3] and Dayanik and Sezer [9], but have never been considered for the combination of Wiener and Poisson processes. Clearly, an unobserved change can be detected more accurately if there are multiple independent sources of information about the disorder time. If all of the information sources consist of exclusively either Wiener or Poisson process observations, then the problem can be solved by applying the results of Shiryaev ([17], Chapter 4) in the Wiener case and Dayanik and Sezer [9] in the Poisson case to a weighted linear combination or superposition of all observation processes; see Section 5. If Wiener and Poisson processes can be observed simultaneously, then previous work does not provide an answer; the solution of the problem in this case is the current paper's contribution.

We solve the problem in detail for $m = d = 1$, namely, when we observe exactly one Wiener and one Poisson process simultaneously; in Section 5 we show the easy extension to multiple Wiener and multiple Poisson processes. Therefore, except in Section 5, we drop all of the superscripts in the sequel. We show that the first time $\tau_{[\phi_\infty,\infty)} \triangleq \inf\{t \geq 0; \Phi_t \geq \phi_\infty\}$ that the conditional odds-ratio process

$$\Phi_t \triangleq \frac{\mathbb{P}\{\Theta \leq t \mid \mathcal{F}_t\}}{\mathbb{P}\{\Theta > t \mid \mathcal{F}_t\}}, \qquad t \geq 0 \tag{1.3}$$

enters into some half-line $[\phi_\infty, \infty) \subset \mathbb{R}_+$ gives the smallest Bayes risk. To calculate the critical threshold $\phi_\infty$ and the minimum Bayes risk, we reduce the original problem to an optimal stopping problem for the process $\Phi$, which turns out to be a jump-diffusion jointly driven by the Wiener and point processes; see (2.8) for its dynamics. The value function of the optimal stopping problem satisfies certain variational inequalities, but they involve a difficult second order integro-differential equation.

We overcome the anticipated difficulties of directly solving the variational inequalities by introducing a jump operator. By means of that operator, we transform the original optimal stopping problem for the jump-diffusion process $\Phi$ into a sequence of optimal stopping problems for the diffusion part $Y$ of the process $\Phi$ between its successive jumps. This decomposition allows us to employ the powerful optimal stopping theory for one-dimensional diffusions to solve each sub-problem between jumps. The solutions of those sub-problems are then combined by means of the jump operator, whose role is basically to incorporate new information about disorder time arriving at jump times of the point process.

Solving optimal stopping problems for jump-diffusion processes by separating jump and diffusion parts with the help of a jump operator seems new and may prove to be useful in other applications, too. Our approach was inspired by several personal conversations with Professor Erhan Çinlar on better ways to calculate the distributions of various functionals of jump processes. For Professor Çinlar's interesting view on this more general matter, his recent lecture [8] in honor of the 2006 Blackwell–Tapia prize recipient, Professor William Massey, in the Blackwell–Tapia conference, held between 3–4 November 2006, may be consulted.

In Section 2 we start our study by giving the precise description of the detection problem and by modeling it under a reference probability measure; the equivalent optimal stopping problem is derived, and the conditional odds-ratio process is examined. In Section 3 we introduce the jump operator. By using it repeatedly, we define "successive approximations" of the optimal stopping problem's value function and identify their important properties. Their common structure is inherited in the limit by the value function and is used at the end of Section 4 to describe an optimal alarm time for the original detection problem. Each successive approximation is itself the value function of some optimal stopping problem, but now for a *diffusion*, and their explicit calculation is undertaken in Section 4. The successive approximations converge uniformly and at an exponential rate to the original value function. Therefore, they are built into an efficient and accurate approximation algorithm, which is explained in Section 6 and illustrated on several examples. Examples suggest that observing Poisson and Wiener processes simultaneously can reduce the Bayes risk significantly. Baron and Tartakovsky [1] have recently derived asymptotic expansions of both optimal critical threshold and minimum Bayes risk as the detection delay cost $c$ tends to zero. In Section 6 we have compared in one of the examples those expansions to the approximations of actual values calculated by our numerical algorithm. Finally, some of the lengthy calculations are deferred to the Appendix.

**2. Problem description and model.** Let $(\Omega, \mathcal{F}, \mathbb{P})$ be a probability space hosting a marked point process $\{(T_n, Z_n); n \geq 1\}$ whose $(E, \mathcal{E})$-valued marks $Z_n$, $n \geq 1$ arrive at times $T_n$, $n \geq 1$, a one-dimensional Wiener process $W$, and a random variable $\Theta$ with distribution in (1.1). The counting measure

$$p((0,t] \times A)) \triangleq \sum_{n=1}^{\infty} \mathbf{1}_{(0,t]\times A}(T_n, Z_n), \qquad t \geq 0, A \in \mathcal{E}$$

generates the internal history $\mathbb{F}^p = \{\mathcal{F}^p\}_{t\geq 0}$,

$$\mathcal{F}_t^p \triangleq \sigma\{p((0,s] \times A); 0 \leq s \leq t, A \in \mathcal{E}\},$$

of the marked point process $\{(T_n, Z_n); n \geq 1\}$. At time $\Theta$, (i) the drift of the Wiener process $W$ changes from zero to $\mu$, and (ii) the $(\mathbb{P}, \mathbb{F}^p)$-compensator of the counting measure $p(dt \times dz)$ changes from $\lambda_0\, dt\nu_0(dz)$ to $\lambda_1\, dt\nu_1(dz)$.

The process $W$ is independent of $\Theta$ and $(T_n, Z_n)_{n\geq 1}$. Neither $W$ nor $\Theta$ are observable. Instead,

$$X_t = X_0 + \mu(t-\Theta)^+ + W_t, \qquad t \geq 0$$

and $\{(T_n, Z_n); n \geq 1\}$ are observable. The observation filtration $\mathbb{F} = \{\mathcal{F}_t\}_{t\geq 0}$ consists of the internal filtrations of $X$ and $(T_n, Z_n)_{n\geq 1}$; that is,

$$\mathcal{F}_t \triangleq \mathcal{F}_t^X \vee \mathcal{F}_t^p \quad \text{and} \quad \mathcal{F}_t^X \triangleq \sigma\{X_s; 0 \leq s \leq t\} \qquad \text{for every } t \geq 0.$$

If we enlarge $\mathbb{F}$ by the information about $\Theta$ and denote the enlarged filtration by $\mathbb{G} = \{\mathcal{G}_t\}_{t\geq 0}$,

$$\mathcal{G}_t \triangleq \mathcal{F}_t \vee \sigma\{\Theta\}, \qquad t \geq 0,$$

then for every nonnegative $\mathbb{G}$-predictable process $\{H(t,z)\}_{t\geq 0}$ indexed by $z \in E$, we have

$$\mathbb{E}\left[\int_{(0,\infty)\times E} H(s,z)p(ds \times dz)\right] = \mathbb{E}\left[\int_0^\infty \int_E H(s,z)\lambda(s,dz)\,ds\right],$$

where $\mathbb{E}$ is the expectation with respect to $\mathbb{P}$, and

$$\lambda(s,dz) \triangleq \lambda_0\nu_0(dz)\mathbf{1}_{[0,\Theta)}(s) + \lambda_1\nu_1(dz)\mathbf{1}_{[\Theta,\infty)}(s), \qquad s \geq 0,$$

is the $(\mathbb{P}, \mathbb{G})$-intensity kernel of the counting measure $p(dt \times dz)$; see Brémaud [5], Chapter VIII.

The rates $0 < \lambda, \lambda_0, \lambda_1 < \infty$, the drift $\mu \in \mathbb{R} \setminus \{0\}$, and the probability measures $\nu_0(\cdot)$, $\nu_1(\cdot)$ on $(E, \mathcal{E})$ are known. The objective is to find a stopping time $\tau$ of the observation filtration $\mathbb{F}$ with the smallest Bayes risk $R_\tau(\pi)$ in (1.2) for every $\pi \in [0, 1)$.

*Model.* Let $(\Omega, \mathcal{F}, \mathbb{P}_0)$ be a probability space hosting the following independent stochastic elements:

(i) a one-dimensional Wiener process $X = \{X_t; t \geq 0\}$,

(ii) an $(E, \mathcal{E})$-valued marked point process $\{(T_n, Z_n); n \geq 1\}$ whose counting measure $p(dt \times dz)$ has $(\mathbb{P}_0, \mathbb{F}^p)$-compensator $\lambda_0\,dt\,\nu_0(dz)$,

(iii) a random variable $\Theta$ with zero-modified exponential distribution

$$\mathbb{P}_0\{\Theta = 0\} = \pi, \qquad \mathbb{P}_0\{\Theta > t\} = (1-\pi)e^{-\lambda t}, \qquad t \geq 0. \tag{2.1}$$

Suppose that $\nu_1(\cdot)$ is absolutely continuous with respect to $\nu_0(\cdot)$ and has the Radon–Nikodym derivative

$$f(z) \triangleq \left.\frac{d\nu_1}{d\nu_0}\right|_{\mathcal{E}}(z), \qquad z \in E. \tag{2.2}$$

Define a new probability measure $\mathbb{P}$ on $\mathcal{G}_\infty = \bigvee_{t\geq 0}\mathcal{G}_t$ locally by means of the Radon–Nikodym derivative of its restriction to $\mathcal{G}_t$,

$$\frac{d\mathbb{P}}{d\mathbb{P}_0}\bigg|_{\mathcal{G}_t} = \xi_t \triangleq \mathbf{1}_{\{t<\Theta\}} + \mathbf{1}_{\{t\geq\Theta\}}\frac{L_t}{L_\Theta}, \qquad t\geq 0,$$

where

$$L_t \triangleq \exp\left\{\mu X_t - \left(\frac{\mu^2}{2}+\lambda_1-\lambda_0\right)t\right\}\prod_{n:0<T_n\leq t}\left(\frac{\lambda_1}{\lambda_0}f(Z_n)\right), \qquad t\geq 0 \tag{2.3}$$

is a likelihood-ratio process with the dynamics $L_0=1$, and

$$dL_t = L_t\mu\,dX_t + L_{t-}\int_E\left(\frac{\lambda_1}{\lambda_0}f(z)-1\right)[p(dt\times dz)-\lambda_0\,dt\,\nu_0(dz)], \qquad t\geq 0. \tag{2.4}$$

Under the probability measure $\mathbb{P}$, the processes $X$ and $\{(T_n,Z_n);n\geq 1\}$ and the random variable $\Theta$ jointly have exactly the same properties as in the above description of the problem. Moreover, the minimum Bayes risk $U(\cdot)$ can be written as

$$U(\pi)\triangleq\inf_{\tau\in\mathbb{F}}R_\tau(\pi) = 1-\pi+c(1-\pi)V\left(\frac{\pi}{1-\pi}\right), \qquad \pi\in[0,1) \tag{2.5}$$

in terms of the value function

$$V(\phi)\triangleq\inf_{\tau\in\mathbb{F}}\mathbb{E}_0^\phi\left[\int_0^\tau e^{-\lambda t}g(\Phi_t)\,dt\right], \qquad \phi\geq 0,\ \text{where } g(\phi)\triangleq\phi-\frac{\lambda}{c} \tag{2.6}$$

of the optimal stopping problem above for the conditional odds-ratio process $\Phi$ in (1.3); see, for example, Bayraktar, Dayanik and Karatzas [2, Proof of Proposition 2.1]. In (2.6), the expectation $\mathbb{E}_0^\phi$ is taken with respect to $\mathbb{P}_0$ conditionally on $\Phi_0=\phi\geq 0$. Bayes formula gives for every $t\geq 0$ that

$$\begin{aligned}\Phi_t &= \frac{\mathbb{E}_0[\xi_t\mathbf{1}_{\{\Theta\leq t\}}\mid\mathcal{F}_t]}{\mathbb{E}_0[\xi_t\mathbf{1}_{\{\Theta>t\}}\mid\mathcal{F}_t]}\\ &= \frac{\mathbb{E}_0[(L_t/L_\Theta)\mathbf{1}_{\{\Theta\leq t\}}\mid\mathcal{F}_t]}{\mathbb{P}_0\{\Theta>t\}}\\ &= \Phi_0e^{\lambda t}L_t + \int_0^t\lambda e^{\lambda(t-s)}\frac{L_t}{L_s}\,ds;\end{aligned} \tag{2.7}$$

by the chain rule and dynamics in (2.4) of the likelihood-ratio process $L$ we find that

$$\begin{aligned}d\Phi_t &= (\lambda+a\Phi_t)\,dt + \Phi_t\mu\,dX_t\\ &\quad + \Phi_{t-}\int_E\left(\frac{\lambda_1}{\lambda_0}f(z)-1\right)p(dt\times dz), \qquad t\geq 0,\end{aligned} \tag{2.8}$$

where

$$a \triangleq \lambda - \lambda_1 + \lambda_0.$$

Let us define for every $k \geq 0$ that $T_0 = T_0^{(k)} \equiv 0$, and

$$\begin{aligned}
X_u^{(k)} &\triangleq X_{T_k+u} - X_{T_k}, \qquad u \geq 0,\\
(T_\ell^{(k)}, Z_\ell^{(k)}) &\triangleq (T_{k+\ell} - T_k, Z_{k+\ell}), \qquad \ell \geq 1,\\
\mathcal{F}_0^{(k)} &\triangleq \sigma\{(T_n, Z_n), 1 \leq n \leq k\} \vee \sigma\{X_v, 0 \leq v \leq T_k\},\\
\mathcal{F}_u^{(k)} &\triangleq \mathcal{F}_0^{(k)} \vee \sigma\{(T_\ell^{(k)}, Z_\ell^{(k)}); 0 < T_\ell^{(k)} \leq u\} \vee \sigma\{X_v^{(k)}, 0 \leq v \leq u\},\\
&\qquad\qquad u \geq 0,\\
L_u^{(k)} &\triangleq \frac{L_{T_k+u}}{L_{T_k}}\\
&= \exp\left\{\mu X_u^{(k)} - \left(\frac{\mu^2}{2} + \lambda_1 - \lambda_0\right)u\right\} \prod_{\ell : 0 < T_\ell^{(k)} \leq u} \left(\frac{\lambda_1}{\lambda_0} f(Z_\ell^{(k)})\right),\\
&\qquad\qquad u \geq 0.
\end{aligned}$$

Then, as in (2.7), we have

$$\begin{aligned}
\Phi_t &= \Phi_{T_k} e^{\lambda(t-T_k)} L_{t-T_k}^{(k)}\\
&\quad + \int_0^{t-T_k} \lambda e^{\lambda(t-T_k-u)} \frac{L_{t-T_k}^{(k)}}{L_u^{(k)}}\, du, \qquad t \geq T_k, k \geq 0,
\end{aligned} \tag{2.9}$$

and $X^{(k)} = \{X_u^{(k)}, u \geq 0\}$ is a $(\mathbb{P}_0, \{\mathcal{F}_u^{(k)}\}_{u\geq 0})$-Wiener process and $\mathbb{P}_0$-independent of the marked point process $\{(T_\ell^{(k)}, Z_\ell^{(k)}); \ell \geq 1\}$, whose local $(\mathbb{P}_0, \{\mathcal{F}_u^{(k)}\}_{u\geq 0})$-characteristics are $(\lambda_0, \nu_0(\cdot))$ for every $k \geq 0$, since $X^{(0)} \equiv X$, $\{(T_\ell^{(0)}, Z_\ell^{(0)})\}_{\ell \geq 1} \equiv \{(T_n, Z_n)\}_{n\geq 1}$, and $\{\mathcal{F}_t^{(0)}\}_{t\geq 0} \equiv \mathbb{F}$. Thus, the first of two implications of (2.9) is that for every Borel $h : \mathbb{R}_+ \mapsto \mathbb{R}_+$ we have

$$\mathbb{E}_0^{\phi}[h(\Phi_t)\mathbf{1}_{\{t\geq T_k\}} \mid \mathcal{F}_{T_k}] = \mathbf{1}_{\{t \geq T_k\}} \mathbb{E}_0^{\Phi_{T_k}}[h(\Phi_v)]|_{v=t-T_k}, \qquad k \geq 0,$$

which also follows from the strong Markov property of the process $\Phi$ applied at the stopping time $T_k$ of the filtration $\mathbb{F}$. To state the second implication, let us introduce the processes

$$\begin{aligned}
Y_u^{k,y} &= y e^{\lambda u} \exp\left\{\mu X_u^{(k)} - \left(\frac{\mu^2}{2} + \lambda_1 - \lambda_0\right)u\right\}\\
&\quad + \int_0^u \lambda e^{\lambda(u-s)} \exp\left\{\mu(X_u^{(k)} - X_s^{(k)})\right.
\end{aligned} \tag{2.10}$$

$$-\left(\frac{\mu^2}{2}+\lambda_1-\lambda_0\right)(u-s)\Bigg\}\,ds$$

for every $u,k,y\geq 0$, which is (2.9) with $u=t-T_k$ and $y=\Phi_{T_k}$ after all of future jumps are stripped away. Then for every $k\geq 0$, the process $\{Y_u^{k,y}, u\geq 0\}$ is a diffusion on $\mathbb{R}_+$ with the dynamics

$$\begin{aligned} dY_t^{k,y} &= (\lambda + aY_t^{k,y})\,dt + \mu Y_t^{k,y}\,dX_t^{(k)}, \qquad t\geq 0 \quad \text{and} \\ Y_0^{k,y} &= y\geq 0, \qquad k\geq 0. \end{aligned} \tag{2.11}$$

Since $X^{(0)}\equiv X$, we shall drop the superscript 0 from $Y^{0,y}$ and denote that process by $Y^y$. Because for every $k\geq 0$, $X^{(k)}$ is a Wiener process, the processes $Y^{k,y}$, $k\geq 0$ and $Y^y$ have the same finite-dimensional distributions, and (2.9) implies that

$$\Phi_t = \begin{cases} Y_{t-T_k}^{k,\Phi_{T_k}}, & t\in[T_k,T_{k+1}),\\ \dfrac{\lambda_1}{\lambda_0} f(Z_{k+1})\Phi_{t-}, & t=T_{k+1}. \end{cases} \tag{2.12}$$

The superscript $k$ in $Y_{t-T_k}^{k,\Phi_{T_k}} \equiv Y_u^{k,y}|_{y=\Phi_{T_k},u=t-T_k}$ indicates that, when $\{Y_u^{k,y}; u\geq 0\}$ is calculated according to (2.10) or (2.11), the increments of the driving Wiener process $\{X_u^{(k)}; u\geq 0\}$ are those of the process $X$ after the $k$th arrival time $T_k$ of the marked point process.

**3. Jump operator.** Let us now go back to the optimal stopping problem in (2.6). By (2.12), we have $\Phi_t = Y_t^{0,\Phi_0}\equiv Y_t^{\Phi_0}$ for $0\leq t<T_1$. This suggests that every $\mathbb{F}$-stopping time $\tau$ coincides on the event $\{\tau<T_1\}$ with one of the stopping times of the process $Y^{\Phi_0}$. On the other hand, the process $\Phi$ regenerates at time $T_1$ starting from its new position $\Phi_{T_1} = (\lambda_1/\lambda_0)f(Z_1)Y_{T_1}^{\Phi_0}$. Moreover, on the event $\{\tau>T_1\}$, the expected total running cost $\mathbb{E}_0^{\phi}[\int_0^{T_1} e^{-\lambda t} g(Y_t^{\Phi_0})\,dt]$ incurred until time $T_1$ is sunken at time $T_1$, and the smallest Bayes risk achievable in the future should be $V(\Phi_{T_1})$ independent of the past. Hence, if we define an operator $J$ acting on the bounded Borel functions $w:\mathbb{R}_+\mapsto\mathbb{R}$ according to

$$(Jw)(\phi) \triangleq \inf_{\tau\in\mathbb{F}^X} \mathbb{E}_0^{\phi}\left[\int_0^{\tau\wedge T_1} e^{-\lambda t} g(\Phi_t)\,dt + \mathbf{1}_{\{\tau\geq T_1\}} e^{-\lambda T_1} w(\Phi_{T_1})\right], \qquad \phi\geq 0, \tag{3.1}$$

then we expect that $V(\phi)=(JV)(\phi)$ for every $\phi\geq 0$. In the next section we prove that $V(\cdot)$ is indeed a fixed point of $w\mapsto Jw$, and if we define $v_n:\mathbb{R}_+\mapsto\mathbb{R}$, $n\geq 0$, successively by

$$v_0(\cdot)\equiv 0 \quad \text{and} \quad v_{n+1}(\cdot)\triangleq (Jv_n)(\cdot), \qquad n\geq 0, \tag{3.2}$$

then $\{v_n(\cdot)\}_{n\geq 1}$ converges to $V(\cdot)$ uniformly. This result will allow us to describe not only an optimal strategy, but also a numerical algorithm that approximates the optimal strategy and the value function.

Note that the infimum in (3.1) is taken over stopping times of the Wiener process $X$. Since $X$ and the marked point process $(T_n, Z_n)_{n\geq 1}$ are $\mathbb{P}_0$-independent, the decomposition in (2.12) and some algebra lead to

$$(Jw)(\phi) = \inf_{\tau\in\mathbb{F}^X} \mathbb{E}_0^\phi\left[\int_0^\tau e^{-(\lambda+\lambda_0)t}(g+\lambda_0(Kw))(Y_t^{\Phi_0})\,dt\right], \qquad \phi\geq 0, \tag{3.3}$$

where $K$ is the operator acting on bounded Borel functions $w:\mathbb{R}_+\mapsto\mathbb{R}$ according to

$$(Kw)(\phi) \triangleq \int_E w\left(\frac{\lambda_1}{\lambda_0}f(z)\phi\right)\nu_0(dz), \qquad \phi\geq 0, \tag{3.4}$$

where $f(\cdot)$ is the Radon–Nikodym derivative in (2.2). The identity in (3.3) shows that $(Jw)(\phi)$ is the value function of an optimal stopping problem for the one-dimensional diffusion $Y^\phi \equiv Y^{0,\phi}$, whose dynamics are given by (2.11). Standard variational arguments imply that, under suitable conditions, the function $Jw(\cdot)$ satisfies

$$\begin{aligned} 0 = \min\{-(Jw)(\phi), [\mathcal{A}_0-(\lambda+\lambda_0)](Jw)(\phi)+g(\phi)+\lambda_0(Kw)(\phi)\},& \\ \phi\geq 0,& \end{aligned} \tag{3.5}$$

where for every twice continuously-differentiable function $\widetilde{w}:\mathbb{R}_+\mapsto\mathbb{R}$,

$$(\mathcal{A}_0\widetilde{w})(\phi) \triangleq \frac{\mu^2}{2}\phi^2\widetilde{w}''(\phi)+(\lambda+a\phi)\widetilde{w}'(\phi) \tag{3.6}$$

is the $(\mathbb{P}_0,\mathbb{F})$-infinitesimal generator of the process $Y^y$, with drift and diffusion coefficients

$$\mu(\phi)\triangleq\lambda+a\phi \quad\text{and}\quad \sigma(\phi)\triangleq\mu\phi, \tag{3.7}$$

respectively. If both $w$ and $Jw$ are replaced with $V$, then (3.5) becomes

$$0=\min\{-V(\phi),(\mathcal{A}-\lambda)V(\phi)+g(\phi)\}, \qquad \phi\geq 0, \tag{3.8}$$

where for every twice continuously-differentiable function $\widetilde{w}:\mathbb{R}_+\mapsto\mathbb{R}$,

$$(\mathcal{A}\widetilde{w})(\phi)\triangleq(\mathcal{A}_0\widetilde{w})(\phi)+\lambda_0[(K-1)\widetilde{w}](\phi), \qquad \phi\geq 0 \tag{3.9}$$

is the $(\mathbb{P}_0,\mathbb{F})$-infinitesimal generator of the process $\Phi$ in (1.3)–(2.8). The identity in (3.8) coincides with the variational inequalities satisfied by the function $V(\cdot)$ of (2.6) under suitable conditions. This coincidence is the second motivation for the introduction of the operator $J$ in (3.1) and for the claim that $V=JV$ must hold.

Reversing the arguments gives additional insight about the role of the operator $J$. If one decides to attack first to the variational inequalities in (3.8) for $V(\cdot)$, then she realizes that solving integro-differential equation $(\mathcal{A}-\lambda)V+g=0$ is difficult. Substituting into (3.8) the decomposition in (3.9) of the operator $\mathcal{A}$ due to diffusion and jump parts gives

$$0=\min\{-V(\phi),[\mathcal{A}_0-(\lambda+\lambda_0)]V(\phi)+g(\phi)+\lambda_0(KV)(\phi)\},\qquad \phi\geq 0.$$

Now $[\mathcal{A}_0-(\lambda+\lambda_0)]V(\phi)+g(\phi)+\lambda_0(KV)(\phi)=0$ is a nonhomogeneous second order ordinary differential equation (ODE) with the forcing function $-g-\lambda(KV)$. If one wants to take full advantage of the rich theory for the solutions of second order ODEs, then she only needs to break the cycle by replacing the unknown $V$ in the forcing function with some known function $w$ and call by $Jw$ the solution of the resulting variational inequalities, namely, (3.5). By repeatedly replacing $w$ with $Jw$, one then hopes that $J^n w$ converges to $V$ as $n\to\infty$. As the next remark shows, the jump operator $J$ can be applied repeatedly to bounded functions, since $Jw$ is bounded whenever $w$ is bounded.

Remark 3.1. For every bounded $w:\mathbb{R}_+\mapsto\mathbb{R}$, the function $Jw:\mathbb{R}_+\mapsto\mathbb{R}_-$ is bounded, and

$$-\left(\frac{\lambda}{c}+\lambda_0\|w^-\|\right)\frac{1}{\lambda+\lambda_0}\leq (Jw)(\cdot)\leq 0,$$

where $\|w^-\|$ is the sup-norm of the negative part of $w(\cdot)$. If $w$ is bounded and $w(\cdot)\geq -1/c$, then $0\geq (Jw)(\cdot)\geq -1/c$. If $w:\mathbb{R}_+\mapsto\mathbb{R}$ is concave, then so is $Jw:\mathbb{R}_+\mapsto\mathbb{R}_-$. The mapping $w\mapsto Jw$ on the collection of bounded functions is monotone.

Proof. Suppose that $w(\cdot)$ is bounded. Since $\tau\equiv 0$ is an $\mathbb{F}^X$-stopping time, we have $(Jw)(\cdot)\leq 0$. Since $(Kw)(\cdot)\geq -\|w^-\|$ and $g(\Phi_t)=\Phi_t-(\lambda/c)\geq -\lambda/c$ for every $t\geq 0$, we have

$$\begin{aligned}0\geq Jw(\phi)&\geq \inf_{\tau\in\mathbb{F}^X}\mathbb{E}_0^{\phi}\int_0^{\tau}e^{-(\lambda+\lambda_0)t}\left(-\frac{\lambda}{c}-\lambda_0\|w^-\|\right)dt\\&=-\left(\frac{\lambda}{c}+\lambda_0\|w^-\|\right)\frac{1}{\lambda+\lambda_0}.\end{aligned}$$

If $w(\cdot)\geq -1/c$, then $\|w^-\|\leq 1/c$ and $0\geq (Jw)(\cdot)\geq -1/c$. Suppose that $w(\cdot)$ is concave. The explicit form in (2.10) indicates that $Y_t^y\equiv Y_t^{0,y}$ is an affine function of $Y_0^y\equiv Y_0^{0,y}=y$. Since $g(\cdot)$ is affine and $(Kw)(\cdot)$ is concave, the mapping $y\mapsto g(Y_t^y)+\lambda_0(Kw)(Y_t^y)$ is also concave. Therefore, the integral in (3.3) and its expectation are concave for every $\mathbb{F}^X$-stopping time $\tau$. Because $(Jw)(\cdot)$ is the infimum of concave functions, it is also concave. The monotonicity of $w\mapsto Jw$ is evident. $\square$

REMARK 3.2. For every $\phi \geq 0$, we have $\mathbb{E}_0^{\phi}[\int_0^{\infty} e^{-\lambda t}\Phi_t\,dt] = \infty$,

$$\mathbb{E}_0^{\phi}\left[\int_0^{\infty} e^{-(\lambda+\lambda_0)t}\Phi_t\,dt\right] = \frac{\phi+1}{\lambda_0} - \frac{1}{\lambda+\lambda_0},$$

$$\mathbb{E}_0^{\phi}\left[\int_0^{\infty} e^{-(\lambda+\lambda_0)t}Y_t^{\Phi_0}\,dt\right] = \frac{1}{\lambda_1}\left(\phi + \frac{1}{\lambda+\lambda_0}\right).$$

PROOF. The proof follows from (2.7), (2.10), Fubini's theorem, and $(\mathbb{P}_0, \mathbb{F})$-martingale property of $L = \{L_t; t \geq 0\}$ after noting that

$$\mathbb{E}_0^{\phi}\Phi_t = (1+\phi)e^{\lambda t} - 1$$

and

$$\mathbb{E}_0^{\phi}Y_t^{\Phi_0} = \phi e^{(\lambda-\lambda_1+\lambda_0)t} + \frac{\lambda[e^{(\lambda-\lambda_1+\lambda_0)t} - 1]}{\lambda-\lambda_1+\lambda_0}. \qquad \square$$

REMARK 3.3. The sequence $\{v_n(\cdot)\}_{n\geq 0}$ in (3.2) is decreasing, and the limit $v_\infty(\phi) \triangleq \lim_{n\to\infty} v_n(\phi)$ exists. The functions $\phi \mapsto v_n(\phi)$, $0 \leq n \leq \infty$, are concave, nondecreasing and bounded between $-1/c$ and zero.

PROOF. We have $v_1(\phi) = (Jv_0)(\phi) \leq 0 \equiv v_0$, since stopping immediately is always possible. Suppose now that $v_n \leq v_{n-1}$ for some $n \geq 1$. Then $v_{n+1} = Jv_n \leq Jv_{n-1} = v_n$ by Remark 3.1, and $\{v_n(\cdot)\}_{n\geq 1}$ is a decreasing sequence by induction. Since $v_0 \equiv 0$ is concave and bounded between 0 and $-1/c$, Remark 3.1 and another induction imply that every $v_n(\cdot)$, $1 \leq n \leq \infty$, is concave and bounded between $-1/c$ and 0. Finally, every concave bounded function on $\mathbb{R}_+$ must be nondecreasing; otherwise, the negative right-derivative at some point does not increase on the right of that point, and the function eventually diverges to $-\infty$. $\square$

LEMMA 3.1. *The function $v_\infty(\cdot) \triangleq \lim_{n\to\infty} v_n(\cdot)$ is the unique bounded solution of the equation $w(\cdot) = (Jw)(\cdot)$.*

PROOF. Since by Remark 3.3 $\{v_n(\cdot)\}_{n\geq 0}$ is a decreasing sequence of bounded functions, the dominated convergence theorem implies that

$$\begin{aligned}
v_\infty(\phi) &= \inf_{n\geq 1} v_{n+1}(\phi) \\
&= \inf_{n\geq 1}\inf_{\tau\in\mathbb{F}^X} \mathbb{E}_0^{\phi}\left[\int_0^{\tau} e^{-(\lambda+\lambda_0)t}[g(Y_t^{\Phi_0}) + \lambda_0(Kv_n)(Y_t^{\Phi_0})]\,dt\right] \\
&= \inf_{\tau\in\mathbb{F}^X}\inf_{n\geq 1} \mathbb{E}_0^{\phi}\left[\int_0^{\tau} e^{-(\lambda+\lambda_0)t}[g(Y_t^{\Phi_0}) + \lambda_0(Kv_n)(Y_t^{\Phi_0})]\,dt\right] \\
&= \inf_{\tau\in\mathbb{F}^X} \mathbb{E}_0^{\phi}\left[\int_0^{\tau} e^{-(\lambda+\lambda_0)t}\left[g(Y_t^{\Phi_0}) + \lambda_0\left(K\inf_{n\geq 1} v_n\right)(Y_t^{\Phi_0})\right]dt\right] \\
&= (Jv_\infty)(\phi).
\end{aligned}$$

Let $u_1(\cdot)$ and $u_2(\cdot)$ be two bounded solutions of $w = Jw$. Fix any arbitrary $\phi \in \mathbb{R}_+$ and $\varepsilon > 0$. Because $(Ju_1)(\phi)$ is finite, there is some $\tau_1 = \tau_1(\phi) \in \mathbb{F}^X$ such that

$$u_1(\phi) = (Ju_1)(\phi) \geq \mathbb{E}_0^{\phi}\left[\int_0^{\tau_1} e^{-(\lambda+\lambda_0)t}(g + \lambda_0(Ku_1))(Y_t^{\Phi_0})\,dt\right] - \varepsilon.$$

Because $Ku_1 - Ku_2 = K(u_1 - u_2) \leq \|u_1 - u_2\|$, we have

$$\begin{aligned}
u_2(\phi) - u_1(\phi) &\leq \mathbb{E}_0^{\phi}\left[\int_0^{\tau_1} e^{-(\lambda+\lambda_0)t}(g + \lambda_0(Ku_2))(Y_t^{\Phi_0})\,dt\right] \\
&\quad - \mathbb{E}_0^{\phi}\left[\int_0^{\tau_1} e^{-(\lambda+\lambda_0)t}(g + \lambda_0(Ku_1))(Y_t^{\Phi_0})\,dt\right] + \varepsilon \\
&= \mathbb{E}_0^{\phi}\left[\int_0^{\tau_1} e^{-(\lambda+\lambda_0)t}\lambda_0(K(u_2 - u_1))(Y_t^{\Phi_0})\,dt\right] + \varepsilon \\
&\leq \|u_2 - u_1\| \int_0^{\infty} \lambda_0 e^{-(\lambda+\lambda_0)t} + \varepsilon \leq \frac{\lambda_0}{\lambda + \lambda_0}\|u_2 - u_1\| + \varepsilon.
\end{aligned}$$

Since $\varepsilon$ is arbitrary, this implies $u_2(\phi) - u_1(\phi) \leq [\lambda_0/(\lambda + \lambda_0)]\|u_2 - u_1\|$. Interchanging $u_1$ and $u_2$ gives $u_1(\phi) - u_2(\phi) \leq [\lambda_0/(\lambda + \lambda_0)]\|u_2 - u_1\|$, and the last two inequalities yield $|u_1(\phi) - u_2(\phi)| \leq [\lambda_0/(\lambda + \lambda_0)]\|u_1 - u_2\|$ for every $\phi \geq 0$. Therefore, $\|u_1 - u_2\| \leq [\lambda_0/(\lambda + \lambda_0)]\|u_1 - u_2\|$, and because $0 < \lambda_0/(\lambda + \lambda_0) < 1$, this is possible if and only if $\|u_1 - u_2\| = 0$; hence, $u_1 \equiv u_2$. Therefore, $w = v_\infty$ is the unique bounded solution of $w = Jw$. □

LEMMA 3.2. *The sequence $\{v_n(\phi)\}_{n\geq 0}$ converges to $v_\infty(\phi)$ as $n \to \infty$ uniformly in $\phi \geq 0$. More precisely, we have*

$$v_\infty(\phi) \leq v_n(\phi) \leq v_\infty(\phi) + \frac{1}{c}\left(\frac{\lambda_0}{\lambda + \lambda_0}\right)^n \qquad \forall n \geq 0, \forall \phi \geq 0. \tag{3.10}$$

PROOF. The first inequality follows from Remark 3.3. We shall prove the second inequality by induction on $n \geq 0$. This inequality is immediate for $n = 0$ since $-1/c \leq v_\infty(\cdot) \leq 0$. Suppose that it is true for some $n \geq 0$. Then induction hypothesis implies that

$$\begin{aligned}
v_{n+1}(\phi) &= \inf_{\tau \in \mathbb{F}^X} \mathbb{E}_0^{\phi}\left[\int_0^{\tau} e^{-(\lambda+\lambda_0)t}[g(Y_t^{\Phi_0}) + \lambda_0(Kv_n)(Y_t^{\Phi_0})]\,dt\right] \\
&\leq \inf_{\tau \in \mathbb{F}^X} \mathbb{E}_0^{\phi}\left[\int_0^{\tau} e^{-(\lambda+\lambda_0)t}\left[g(Y_t^{\Phi_0})\right.\right. \\
&\qquad\qquad \left.\left. + \lambda_0(Kv_\infty)(Y_t^{\Phi_0}) + \frac{\lambda_0}{c}\left(\frac{\lambda_0}{\lambda + \lambda_0}\right)^n\right]dt\right] \\
&\leq \inf_{\tau \in \mathbb{F}^X}\left(\mathbb{E}_0^{\phi}\left[\int_0^{\tau} e^{-(\lambda+\lambda_0)t}[g(Y_t^{\Phi_0}) + \lambda_0(Kv_\infty)(Y_t^{\Phi_0})]\,dt\right]\right.
\end{aligned}$$

$$+ \int_0^\infty e^{-(\lambda+\lambda_0)t} \frac{\lambda_0}{c} \left( \frac{\lambda_0}{\lambda+\lambda_0} \right)^n dt \Bigg)$$

$$= (Jv_\infty)(\phi) + \int_0^\infty e^{-(\lambda+\lambda_0)t} \frac{\lambda_0}{c} \left( \frac{\lambda_0}{\lambda+\lambda_0} \right)^n dt$$

$$= v_\infty(\phi) + \frac{1}{c} \left( \frac{\lambda_0}{\lambda+\lambda_0} \right)^{n+1},$$

since $v_\infty = Jv_\infty$ by Lemma 3.1. □

**4. Solution of the optimal stopping problem.** The main results of this section are that $v_\infty(\cdot)$ coincides with the value function $V(\cdot)$ of the optimal stopping problem in (2.6), and that the first entrance time of the process $\Phi$ of (1.3) into half line $[\phi_\infty, \infty)$ for some constant $\phi_\infty > 0$ is optimal for (2.6). We also describe $\varepsilon$-optimal $\mathbb{F}$-stopping times for (2.6) and summarize the calculation of its value function $V(\cdot)$.

We shall first find an explicit solution of the optimal stopping problem in (3.3). The second order ODE $(\lambda + \lambda_0)h(\cdot) = \mathcal{A}_0 h(\cdot)$ on $(0, \infty)$ admits two twice-continuously differentiable solutions, $\psi(\cdot)$ and $\eta(\cdot)$, unique up to multiplication by a positive constant, such that they are increasing and decreasing, respectively. For this and other facts below about one-dimensional diffusions, see, for example, Itô and McKean [11], Borodin and Salminen [4] Karlin and Taylor [12], Chapter 15.

The explicit form in (2.10) of the process $Y^y \equiv Y^{0,y}$ suggests that the process may start at $y = 0$, but then moves instantaneously into $(0, \infty)$ without ever coming back to 0. It can neither start at nor reach from inside to the right boundary located at $\infty$. Indeed, calculated in terms of the *scale function* $S(\cdot)$ and *speed measure* $M(\cdot)$, defined respectively by

$$S(dy) \triangleq \exp\left\{ -2 \int_c^y \frac{\mu(u)}{\sigma^2(u)} du \right\} dy, \qquad y > 0 \quad \text{and}$$
(4.1)
$$M(dy) \triangleq \frac{dy}{\sigma^2(y) S'(y)}, \qquad y > 0$$

for some arbitrary but fixed constant $c > 0$, Feller's boundary tests give

$$S(0+) = -\infty \quad \text{and} \quad \int_0^c \int_0^z M(dy) S(dz) < \infty, \tag{4.2}$$

$$\int_c^\infty \int_z^\infty S(dy) M(dz) = \infty \qquad \text{and} \qquad \int_c^\infty \int_z^\infty M(dy) S(dz) = \infty, \tag{4.3}$$

as shown in Appendix A.1, and according to Table 6.2 of Karlin and Taylor ([12], page 234), we conclude that $y = 0$ and $y = \infty$ are *entry-not-exit* and *natural* boundaries of the state-space $[0,\infty)$, respectively. Therefore, $\psi(\cdot)$ and $\eta(\cdot)$ satisfy boundary conditions

$$\begin{aligned}
&0<\psi(0+)<\infty, && \eta(0+)=\infty,\\
&\lim_{y\to 0+}\frac{\psi'(y)}{S'(y)}=0, && \lim_{y\to 0+}\frac{\eta'(y)}{S'(y)}>-\infty,\\
&\psi(\infty)=\infty, && \eta(\infty)=0,\\
&\lim_{y\to\infty}\frac{\psi'(y)}{S'(y)}=\infty, && \lim_{y\to\infty}\frac{\eta'(y)}{S'(y)}=0.
\end{aligned} \tag{4.4}$$

We shall set $\psi(0)=\psi(0+)$ and $\eta(0)=\eta(0+)$. The Wronskian $B(\cdot)$ of $\psi(\cdot)$ and $\eta(\cdot)$ equals

$$B(y)\triangleq\psi'(y)\eta(y)-\psi(y)\eta'(y)=B(c)S'(y),\qquad y>0, \tag{4.5}$$

where the constant $c$ and that in the scale function $S(\cdot)$ in (4.1) are the same. The second equality is obtained by solving the differential equation $\mathcal{A}_0B=0$, which follows from the equations $\mathcal{A}_0\psi=(\lambda+\lambda_0)\psi$ and $\mathcal{A}_0\eta=(\lambda+\lambda_0)\eta$ after first multiplying these respectively with $\eta$ and $\psi$, and then, subtracting from each other. Observe that

$$B(c)=\frac{B(y)}{S'(y)}=\frac{\psi'(y)}{S'(y)}\eta(y)-\psi(y)\frac{\eta'(y)}{S'(y)},\qquad y\geq 0$$

is constant. Dividing (4.5) by $-\psi^2(y)$ and then integrating the equation give

$$\frac{\eta(y)}{\psi(y)}=\frac{\eta(c)}{\psi(c)}-\int_c^y B(c)\frac{S'(z)}{\psi^2(z)}\,dz,\qquad y\geq 0. \tag{4.6}$$

This identity implies that the constant $B(c)$ must be strictly positive, since the functions $\psi(\cdot)$ and $\eta(\cdot)$ are linearly independent [note that their nontrivial linear combinations cannot vanish at 0 because of (4.4)].

For every Borel subset $D$ of $\mathbb{R}_+$, denote the first entrance time of $Y^y$ and $\Phi$ to $D$ by

$$\tau_D\triangleq\inf\{t\geq 0: Y_t^y\in D\}\quad\text{and}\quad\widetilde{\tau}_D\triangleq\inf\{t\geq 0:\Phi_t\in D\}, \tag{4.7}$$

respectively. If $D=\{z\}$ for some $z\in\mathbb{R}_+$, we will use $\tau_z(\widetilde{\tau}_z)$ instead of $\tau_{\{z\}}(\widetilde{\tau}_{\{z\}})$. Then

$$\begin{aligned}
\mathbb{E}_0^y[e^{-(\lambda+\lambda_0)\tau_z}]&=\frac{\psi(y)}{\psi(z)}\cdot\mathbf{1}_{(0,z]}(y)+\frac{\eta(y)}{\eta(z)}\cdot\mathbf{1}_{(z,\infty)}(y)\\
&\qquad\forall z>0,\forall y\geq 0,
\end{aligned} \tag{4.8}$$

which can be obtained by applying the optional sampling theorem to the $(\mathbb{P}_0, \mathbb{F})$-martingales $\{e^{-(\lambda+\lambda_0)t}\psi(Y_t^y); t \geq 0\}$ and $\{e^{-(\lambda+\lambda_0)t}\eta(Y_t^y); t \geq 0\}$. For every fixed real number $z > 0$, (4.8) implies that

$$\psi(y) = \left\{ \begin{array}{ll} \psi(z)\mathbb{E}_0^y[e^{-(\lambda+\lambda_0)\tau_z}], & 0 \leq y \leq z \\ \dfrac{\psi(z)}{\mathbb{E}_z[e^{-(\lambda+\lambda_0)\tau_y}]}, & y > z \end{array} \right\},$$

$$\eta(y) = \left\{ \begin{array}{ll} \dfrac{\eta(z)}{\mathbb{E}_z[e^{-(\lambda+\lambda_0)\tau_y}]}, & 0 \leq y \leq z \\ \eta(z)\mathbb{E}_0^y[e^{-(\lambda+\lambda_0)\tau_z}], & y > z \end{array} \right\},$$

and suggests a way to calculate functions $\psi(\cdot)$ and $\eta(\cdot)$ up to a multiplication by a constant on a lattice inside $(0, z]$ by using simulation methods. Let us set $\psi(z) = \eta(z) = 1$ (or to any arbitrary positive constant), and suppose that the grid size $h > 0$ and some integer $N$ are chosen such that $Nh = z$. Let $z_n = nh$, $n = 0, \ldots, N$. Then (4.8) implies that one can calculate

$$\begin{aligned} \psi(z_n) &= \psi(z_{n+1})\mathbb{E}_0^{z_n}[\exp\{-(\lambda+\lambda_0)\tau_{z_{n+1}}\}], \qquad n = N-1, \ldots, 1, 0, \\ \eta(z_n) &= \eta(z_{n+1})/\mathbb{E}_0^{z_{n+1}}[\exp\{-(\lambda+\lambda_0)\tau_{z_n}\}], \qquad n = N-1, \ldots, 1, 0, \end{aligned} \tag{4.9}$$

backward from $z_N \equiv z$ by evaluating expectations using simulation.

The functions $\psi(\cdot)$ and $\eta(\cdot)$ can also be characterized as power series or Kummer's functions; see Polyanin and Zaitsev ([16], pages 221, 225, 229, Equation 134 in Section 2.1.2). Those functions take simple forms for certain values of $\lambda$, $\lambda_1$, $\lambda_0$ and $\mu$. For example, if $a = \lambda + \lambda_0 - \lambda_1 \geq 0$ and

$$(n-1)\lambda = (n-2)[(\lambda+\lambda_0) + \tfrac{1}{2}\mu^2(n-1)] \qquad \text{for some } n \in \mathbb{N} \text{ and } n > 2,$$

then $\psi(\cdot)$ is a polynomial of the form $\psi(\phi) = \sum_{k=0}^{n-1} \beta_k \phi^k$, where $\beta_0 = 1$, $\beta_1 = (\lambda+\lambda_0)/\lambda$, and

$$\beta_k = \left[\frac{(\lambda+\lambda_0) - (k-1)a - 0.5\mu^2(k-1)(k-2)}{k\lambda}\right]\beta_{k-1} \qquad \text{for } k \geq 2,$$

and $\eta(\cdot)$ can be obtained in terms of $\psi(\cdot)$ from (4.6). However, we make no such assumptions about the parameters and work with general $\psi(\cdot)$ and $\eta(\cdot)$.

LEMMA 4.1. *Every moment of the first entrance times* $\tau_{[r,\infty)}$ *and* $\widetilde{\tau}_{[r,\infty)}$ *of the processes* $Y^{\Phi_0}$ *and* $\Phi$, *respectively, into half line* $[r, \infty)$ *is uniformly bounded for every* $r \geq 0$.

PROOF. Fix $r > 0$ and $0 \le \phi < r$; the cases $r = 0$ or $\phi \ge r$ are obvious. Since the sample paths of $Y^{\Phi_0}$ are continuous, we have $\tau_{[r,\infty)} \equiv \tau_r$, and (4.8) implies that

$$\begin{aligned} \mathbb{P}_0^\phi\{\tau_r < T_1\} &= \mathbb{E}_0^\phi e^{-\lambda_0 \tau_r} \ge \mathbb{E}_0^\phi e^{-(\lambda+\lambda_0)\tau_r} \\ &= \frac{\psi(\phi)}{\psi(r)} \ge \frac{\psi(0)}{\psi(r)} \in (0,1), \qquad \phi \in [0,r). \end{aligned} \tag{4.10}$$

Let $\alpha \triangleq \sqrt{1 - (\psi(0)/\psi(r))} < 1$. The strong $(\mathbb{P}_0, \mathbb{F})$-Markov property of $Y^{\Phi_0}$ implies that

$$\begin{aligned} &\mathbb{P}_0^\phi\{\tau_r > T_n\} \\ &\quad = \mathbb{P}_0^\phi\{\tau_r > T_{n-1}, \tau_r > T_n\} \\ &\quad = \mathbb{E}_0^\phi[\mathbf{1}_{\{\tau_r > T_{n-1}\}}(\mathbf{1}_{\{\tau_r > T_1\}} \circ \theta_{T_{n-1}})] = \mathbb{E}_0^\phi[\mathbf{1}_{\{\tau_r > T_{n-1}\}} \mathbb{P}_0^{Y^{\Phi_0}_{T_{n-1}}}\{\tau_r > T_1\}] \\ &\quad \le \mathbb{P}_0^\phi\{\tau_r > T_{n-1}\}\left[1 - \frac{\psi(0)}{\psi(r)}\right] \le \left[1 - \frac{\psi(0)}{\psi(r)}\right]^n = \alpha^{2n} \end{aligned} \tag{4.11}$$

by induction on $n$, because $Y^{\Phi_0}_{T_n} \in [0,r)$ on $\{\tau_r > T_n\}$ for every $n \ge 1$. For every $k \ge 1$,

$$\begin{aligned} \mathbb{E}_0^\phi \tau_r^k &\le \mathbb{E}_0^\phi \sum_{n=0}^\infty T_{n+1}^k \mathbf{1}_{\{T_n < \tau_r \le T_{n+1}\}} \\ &\le \sum_{n=0}^\infty \mathbb{E}_0^\phi T_{n+1}^k \mathbf{1}_{\{\tau_r > T_n\}} \le \sum_{n=0}^\infty \sqrt{\mathbb{E}_0^\phi T_{n+1}^{2k} \mathbb{P}_0^\phi\{\tau_r > T_n\}} \\ &\le \lambda_0^{-k} \sum_{n=0}^\infty \sqrt{\frac{(n+2k)!}{n!}} \alpha^n \le \lambda_0^{-k} \sum_{n=0}^\infty (n+2k)^k \alpha^n < \infty \end{aligned} \tag{4.12}$$

independent of the initial state $\phi \ge 0$. Since $\mathbb{P}_0^\phi\{\widetilde{\tau}_{[r,\infty)} < T_1\} = \mathbb{P}_0^\phi\{\tau_r < T_1\} \ge \psi(0)/\psi(r)$ for every $\phi \in [0,r)$ by (4.10), both (4.11) and (4.12) remain correct if we replace $\tau_r$ and $Y^{\Phi_0}_{T_n-1}$ with $\widetilde{\tau}_{[r,\infty)}$ and $\Phi_{T_{n-1}}$, respectively. $\square$

ASSUMPTION. In the remainder, suppose that $w : \mathbb{R}_+ \mapsto \mathbb{R}$ is an arbitrary but fixed bounded and continuous function, and $0 < l < r < \infty$.

Define

$$(H_{l,r}w)(\phi) \triangleq \mathbb{E}_0^\phi\left[\int_0^{\tau_{[0,l]} \wedge \tau_{[r,\infty)}} e^{-(\lambda+\lambda_0)t}(g + \lambda_0(Kw))(Y_t^{\Phi_0})\, dt\right], \qquad \phi \ge 0, \tag{4.13}$$

$$(H_r w)(\phi) \triangleq \mathbb{E}_0^{\phi}\left[\int_0^{\tau_{[r,\infty)}} e^{-(\lambda+\lambda_0)t}(g+\lambda_0(Kw))(Y_t^{\Phi_0})\,dt\right], \qquad \phi \geq 0.$$

We shall first derive the analytical expression below in (4.16) for $(H_r w)(\cdot)$. Since the left boundary at 0 is entrance-not-exit for the process $Y^{\Phi_0}$, that boundary is inaccessible from the interior $(0,\infty)$ of the state-space, and $\lim_{l\searrow 0}\tau_l \wedge \tau_r = \tau_r$ $\mathbb{P}_0^{\phi}$-a.s. for every $\phi > 0$. Because $(Kw)(\cdot)$, is bounded, and $g(\phi) = \phi - (\lambda/c)$, $\phi \geq 0$, is bounded from below, Remark 3.2 and the monotone convergence theorem imply that

$$\begin{aligned}(H_r w)(\phi) &= \lim_{l\searrow 0}(H_{l,r}w)(\phi), \qquad \phi > 0 \quad \text{and}\\ (H_r w)(0) &= \lim_{\phi\searrow 0}\lim_{l\searrow 0}(H_{l,r}w)(\phi)\end{aligned} \tag{4.14}$$

follows from the strong Markov property; see Appendix A.2 for the details. By means of the first equality, the second becomes $(H_r w)(0) = \lim_{\phi\searrow 0}(H_r w)(\phi)$, that is, the function $\phi \mapsto (H_r w)(\phi)$ is continuous at $\phi = 0$. In terms of the fundamental solutions

$$\psi_l(y) \triangleq \psi(y) - \frac{\psi(l)}{\eta(l)}\eta(y) \quad \text{and} \quad \eta_r(y) \triangleq \eta(y) - \frac{\eta(r)}{\psi(r)}\psi(y)$$

of the equation $[\mathcal{A}_0 - (\lambda+\lambda_0)]h(y) = 0$, $l < y < r$ with boundary conditions $h(l) = 0$ and $h(r) = 0$, respectively, and their Wronskian

$$B_{l,r}(y) \triangleq \psi_l'(y)\eta_r(y) - \psi_l(y)\eta_r'(y) = B(y)\left[1 - \frac{\psi(l)}{\eta(l)}\frac{\eta(r)}{\psi(r)}\right],$$

we find, as shown in Appendix A.3, that

$$\begin{aligned}(H_{l,r}w)(\phi) &= \psi_l(\phi)\int_\phi^r \frac{2\eta_r(z)}{\sigma^2(z)B_{l,r}(z)}(g+\lambda_0(Kw))(z)\,dz\\ &\quad + \eta_r(\phi)\int_l^\phi \frac{2\psi_l(z)}{\sigma^2(z)B_{l,r}(z)}(g+\lambda_0(Kw))(z)\,dz, \qquad 0 < l \leq \phi \leq r,\end{aligned} \tag{4.15}$$

where $\sigma(z) = \mu z$ is the diffusion coefficient of the process $Y^{\Phi_0}$ in (3.7). After taking the limit as $l \searrow 0$, the monotone convergence and boundary conditions in (4.4) give

$$\begin{aligned}(H_r w)(\phi) &= \psi(\phi)\int_\phi^r \frac{2\eta(z)}{\sigma^2(z)B(z)}(g+\lambda_0(Kw))(z)\,dz\end{aligned} \tag{4.16}$$

$$+\,\eta(\phi)\int_0^\phi \frac{2\psi(z)}{\sigma^2(z)B(z)}(g+\lambda_0(Kw))(z)\,dz$$
$$-\,\psi(\phi)\frac{\eta(r)}{\psi(r)}\int_0^r \frac{2\psi(z)}{\sigma^2(z)B(z)}(g+\lambda_0(Kw))(z)\,dz, \qquad 0<\phi\le r,$$

and $(H_r w)(0) = \lim_{\phi\searrow 0}(H_r w)(\phi)$ by (4.14). Finally, $(H_r w)(\phi) = 0$ for every $\phi > r$ by the definition in (4.13). For every $r > 0$, the function $\phi \mapsto (H_r w)(\phi)$ is continuous on $[0,\infty)$; it is twice continuously-differentiable on $(0,\infty)$, possibly except at $\phi = r$. Direct calculation shows that $(H_r w)(r) = (H_r w)'(r+) = 0$ and

$$(H_r w)'(r-) = \left[\eta'(r) - \frac{\eta(r)}{\psi(r)}\psi'(r)\right]\int_0^r \frac{2\psi(z)}{\sigma^2(z)B(z)}(g+\lambda_0(Kw))(z)\,dz.$$

Since $z \mapsto \eta(z) - [\eta(r)/\psi(r)]\psi(z)$ is strictly decreasing,

$$(4.17)\qquad \begin{aligned}&(H_r w)'(r-) = 0\\ &\iff \quad (Gw)(r) \triangleq \int_0^r \frac{2\psi(z)}{\sigma^2(z)B(z)}(g+\lambda_0(Kw))(z)\,dz = 0.\end{aligned}$$

Lemma 4.2. *If $w(\cdot)$ is nondecreasing and nonpositive, then $(Gw)(\phi) = 0$ has exactly one strictly positive solution $\phi = \phi[w]$. If we denote by $\phi_\ell[w]$ the unique solution $\phi$ of $(g + \lambda_0(Kw))(\phi) = 0$ and define $\phi_r[w] \triangleq \phi[-\|w\|]$, then $\phi_\ell[w] \le \phi[w] \le \phi_r[w]$. Moreover, $(Gw)(\phi)$ is strictly negative for $\phi \in (0, \phi[w])$ and strictly positive for $\phi \in (\phi[w], \infty)$.*

Proof. Since $\phi \mapsto (g+\lambda_0(Kw))(\phi) = \phi - (\lambda/c) + \lambda_0(Kw)(\phi)$ is negative at $\phi = 0$ and increases unboundedly as $\phi \to \infty$, it has unique root at some $\phi = \phi_\ell[w] > 0$. Therefore,

$$(Gw)'(\phi) = \frac{2\psi(\phi)}{\sigma^2(\phi)B(\phi)}(g+\lambda_0(Kw))(\phi)$$

changes its sign exactly once at $\phi = \phi_\ell[w]$, from negative to positive, and the continuously differentiable function $(Gw)(\phi) = \int_0^\phi (Gw)'(z)\,dz$ is strictly negative on $(0, \phi_\ell[w]]$. Since $(Gw)(\phi)$ is increasing at every $\phi \in [\phi_\ell[w], \infty)$, the proof will be complete if we show that $\lim_{\phi\to\infty}(Gw)'(\phi) = \infty$. Since $\sigma^2(\phi) = \mu^2\phi^2$, and

$$S'(\phi) = \exp\left\{-2\int_c^\phi \frac{\lambda + au}{\mu^2 u^2}\,du\right\} = \text{const.} \times \exp\left\{\frac{2\lambda}{\mu^2\phi}\right\}\phi^{-2a/\mu^2},$$

we have

$$\lim_{\phi\to\infty}(Gw)'(\phi) = \lim_{\phi\to\infty}\frac{2\psi(\phi)\phi}{\sigma^2(\phi)B(\phi)} = \text{const.} \times \lim_{\phi\to\infty}\frac{\psi(\phi)}{\phi^{1-(2a/\mu^2)}},$$

which equals $\infty$ if $1-(2a/\mu^2)\le 0$. Otherwise, the L'Hospital rule and (4.4) give

$$\lim_{\phi\to\infty}(Gw)'(\phi)=\text{const.}\times\lim_{\phi\to\infty}\frac{\psi'(\phi)}{\phi^{-2a/\mu^2}}$$
$$=\text{const.}\times\lim_{\phi\to\infty}\frac{\psi'(\phi)}{S'(\phi)}\exp\left\{-\frac{2\lambda}{\mu^2\phi}\right\}=\infty.$$

Finally, constant function $w_0(\phi)\triangleq-\|w\|$, $\phi\ge 0$, is also bounded continuous nondecreasing and nonpositive. By the first part of the lemma, $(Gw_0)(\phi)=0$ has exactly one strictly positive solution $\phi=\phi[w_0]=:\phi_r[w]$. Since $w(\cdot)\ge-\|w\|$, we have $(Gw)(\cdot)\ge(Gw_0)(\cdot)$, and therefore, $\phi[w]\le\phi_r[w]$. □

Lemma 4.2 and (4.17) show that in the family of functions $\{H_r(\phi),\phi\in\mathbb{R}_+\}_{r>0}$ there is exactly one function that "fits smoothly at $\phi=r$" and is therefore continuously differentiable on the whole $\phi\in(0,\infty)$, and that function corresponds to the unique strictly positive solution $r=\phi[w]$ of the equation $(Gw)(r)=0$ in (4.17).

Lemma 4.3. *Suppose that $w(\cdot)$ is nondecreasing and nonpositive. Then the function*

$$(Hw)(\phi)\triangleq(H_{\phi[w]}w)(\phi),\qquad \phi\ge 0,\tag{4.18}$$

*equals zero for $\phi>\phi[w]$ and*

$$\psi(\phi)\int_\phi^{\phi[w]}\frac{2\eta(z)}{\sigma^2(z)B(z)}(g+\lambda_0(Kw))(z)\,dz$$
$$+\eta(\phi)\int_0^\phi\frac{2\psi(z)}{\sigma^2(z)B(z)}(g+\lambda_0(Kw))(z)\,dz$$

*for $0<\phi\le\phi[w]$. It is bounded continuous on $[0,\infty)$, continuously differentiable on $(0,\infty)$ and twice continuously differentiable on $(0,\infty)\setminus\{\phi[w]\}$. It satisfies $(Hw)(\phi[w])=(Hw)'(\phi[w])=0$ and the variational inequalities*

$$\left\{\begin{array}{r}(Hw)(\phi)<0\\ [\mathcal{A}_0-(\lambda+\lambda_0))](Hw)(\phi)+(g+\lambda_0(Kw))(\phi)=0\end{array}\right\},\qquad \phi\in(0,\phi[w]),\tag{4.19}$$

$$\left\{\begin{array}{r}(Hw)(\phi)=0\\ [\mathcal{A}_0-(\lambda+\lambda_0)](Hw)(\phi)+(g+\lambda_0(Kw))(\phi)>0\end{array}\right\},\qquad \phi\in(\phi[w],\infty).\tag{4.20}$$

PROOF. The explicit form of $(Hw)(\cdot)$ follows from (4.16) after noticing that the third term equals $-\psi(\phi)[\eta(r)/\psi(r)](Gw)(r)$ and vanishes for $r=\phi[w]$ by definition. Since $(H_r w)(\cdot)$ is continuous on $[0,\infty)$ and twice continuously differentiable on $(0,\infty)\setminus\{r\}$ and $(H_r w)(r)=0$ for every $r>0$, so is $(Hw)(\cdot)\equiv(H_{\phi[w]}w)(\cdot)$ and $(Hw)(\phi[w])=0$. It is also continuously differentiable at $\phi=\phi[w]$ since $(Hw)'(\phi[w]-)\equiv(H_{\phi[w]}w)'(\phi[w]-)=0=(Hw)'(\phi[w]+)$ by (4.17) and Lemma 4.2. Because the function $(Hw)(\cdot)$ is continuous everywhere and vanishes outside the closed and bounded interval $[0,\phi[w]]$, it is bounded everywhere. Direct calculation gives immediately the equalities in (4.19) and (4.20). The inequality in (4.20) follows from substitution of $(Hw)(\phi)=0$ for $\phi>\phi[w]$ and that $(g+\lambda_0(Kw))(\phi)>0$ for $\phi>\phi[w]>\phi_\ell[w]$ by Lemma 4.2, where $\phi_\ell[w]$ is the unique root of nondecreasing function $\phi\mapsto(g+\lambda_0(Kw))(\phi)$. For the proof of the inequality in (4.19), note that $(Hw)'(\phi)$ equals

$$\psi'(\phi)\int_\phi^{\phi[w]}\frac{2\eta(z)}{\sigma^2(z)B(z)}(g+\lambda_0(Kw))(z)\,dz \\ +\eta'(\phi)\int_0^\phi\frac{2\psi(z)}{\sigma^2(z)B(z)}(g+\lambda_0(Kw))(z)\,dz$$

for $0<\phi\le\phi[w]$. The second term is positive since (i) $\eta(\cdot)$ is strictly decreasing, and (ii) $(Gw)(\phi)$ in (4.17) is strictly negative for $\phi\in(0,\phi[w])$ by Lemma 4.2. The first term is strictly negative for $\phi\in(\phi_\ell[w],\phi[w])$, since (i) $\psi(\cdot)$ is strictly increasing, and (ii) $(g+\lambda_0(Kw))(z)>0$ for $z>\phi_\ell[w]$. Therefore, $(Hw)'(\phi)>0$ for $\phi\in[\phi_\ell[w],\phi[w])$. Because continuously differentiable $(Hw)(\phi)$ vanishes at $\phi=\phi[w]$, we have

$$(Hw)(\phi)=-\int_\phi^{\phi[w]}(Hw)'(z)\,dz<0 \qquad \text{for every } \phi_\ell[w]\le\phi<\phi[w].$$

Finally, for every $0\le\phi\le\phi_\ell[w]$, the strong Markov property of the process $Y^{\Phi_0}$ applied at the $\mathbb{F}$-stopping time $\tau_{\phi_\ell[w]}$ gives

$$(Hw)(\phi)=\mathbb{E}_0^\phi\left[\int_0^{\tau_{\phi_\ell[w]}}e^{-(\lambda+\lambda_0)t}(g+\lambda_0(Kw))(Y_t^{\Phi_0})\,dt\right] \\ +\mathbb{E}_0^\phi[e^{-(\lambda+\lambda_0)\tau_{\phi_\ell[w]}}](Hw)(\phi_\ell[w]),$$

and both terms are strictly negative, since $(g+\lambda_0(Kw))(\phi)<0$ for $\phi\in[0,\phi_\ell[w])$ and $(Hw)(\phi_\ell[w])<0$ by the previous displayed equation. □

PROPOSITION 4.1. *Suppose $w(\cdot)$ is nondecreasing and nonpositive. Then*

$$(Jw)(\phi)=(Hw)(\phi)\equiv\mathbb{E}_0^\phi\left[\int_0^{\tau_{[\phi[w],\infty)}}e^{-(\lambda+\lambda_0)t}(g+\lambda_0(Kw))(Y_t^{\Phi_0})\,dt\right], \\ \phi\ge 0.$$

PROOF. For every $0<l<\phi<r$ and $\mathbb{F}^X$-stopping time $\tau$, Itô's rule yields

$$\begin{aligned}
&e^{-(\lambda+\lambda_0)(\tau\wedge\tau_l\wedge\tau_r)}(Hw)(Y^{\Phi_0}_{\tau\wedge\tau_l\wedge\tau_r})\\
&\quad=(Hw)(\phi)\\
&\qquad+\int_0^{\tau\wedge\tau_l\wedge\tau_r} e^{-(\lambda+\lambda_0)t}\mu Y_t^{\Phi_0}(Hw)'(Y_t^{\Phi_0})\,dX_t\\
&\qquad+\int_0^{\tau\wedge\tau_l\wedge\tau_r} e^{-(\lambda+\lambda_0)t}[\mathcal{A}_0-(\lambda+\lambda_0)](Hw)(Y_t^{\Phi_0})\,dt.
\end{aligned}$$

Since $(Hw)'(\cdot)$ is continuous by Lemma 4.3, it is bounded on $[l,r]$. Taking expectations gives

$$\begin{aligned}
&\mathbb{E}_0^{\phi}[e^{-(\lambda+\lambda_0)(\tau\wedge\tau_l\wedge\tau_r)}(Hw)(Y^{\Phi_0}_{\tau\wedge\tau_l\wedge\tau_r})]\\
&\quad=(Hw)(\phi)+\mathbb{E}_0^{\phi}\left[\int_0^{\tau\wedge\tau_l\wedge\tau_r} e^{-(\lambda+\lambda_0)t}[\mathcal{A}_0-(\lambda+\lambda_0)](Hw)(Y_t^{\Phi_0})\,dt\right]\\
&\quad\ge(Hw)(\phi)-\mathbb{E}_0^{\phi}\left[\int_0^{\tau\wedge\tau_l\wedge\tau_r} e^{-(\lambda+\lambda_0)t}(g+\lambda_0(Kw))(Y_t^{\Phi_0})\,dt\right],
\end{aligned}$$

because $(Hw)(\cdot)$ satisfies the variational inequalities in (4.19) and (4.20) by Lemma 4.3. Since $(Hw)(\cdot)\equiv(H_{\phi[w]})(\cdot)$ is nonpositive continuous and bounded by the same lemma, letting $l\to0$, $r\to\infty$ and the dominated convergence theorem (see Remark 3.2) give

$$\begin{aligned}
0&\ge\mathbb{E}_0^{\phi}[e^{-(\lambda+\lambda_0)\tau}(Hw)(Y_\tau^{\Phi_0})]\\
&\ge(Hw)(\phi)-\mathbb{E}_0^{\phi}\left[\int_0^{\tau} e^{-(\lambda+\lambda_0)t}(g+\lambda_0(Kw))(Y_t^{\Phi_0})\,dt\right].
\end{aligned}$$

Thus, we have

$$\mathbb{E}_0^{\phi}\left[\int_0^{\tau} e^{-(\lambda+\lambda_0)t}(g+\lambda_0(Kw))(Y_t^{\Phi_0})\,dt\right]\ge(Hw)(\phi).$$

Taking infimum over $\mathbb{F}^X$-stopping times $\tau$ gives $(Jw)(\phi)\ge(Hw)(\phi)$, $\phi>0$.

If we replace every $\tau$ above with the first entrance time $\tau_{[\phi[w],\infty)}$ of the process $Y^{\Phi_0}$ into $[\phi[w],\infty)$, then $\mathbb{P}_0^{\phi}\{\tau<\infty\}=1$ and the variational inequalities in (4.19) and (4.20) ensure that every inequality above becomes an equality. This proves $(Jw)(\phi)=(Hw)(\phi)$ for every $\phi>0$. Finally, that equality extends to $\phi=0$ by the continuity of $(Jw)(\cdot)$ and $(Hw)(\cdot)$. □

COROLLARY 4.1. *Recall the sequence $\{v_n(\cdot)\}_{n\ge0}$ of functions defined successively by (3.2) and its pointwise limit $v_\infty(\cdot)$, all of which are bounded, concave, nonpositive and nondecreasing by Remark 3.3. Then every $v_n(\cdot)$,*

*$0 \le n \le \infty$, is continuous on $[0,\infty)$ continuously differentiable on $(0,\infty)$, and twice continuously differentiable on $(0,\infty)\setminus\{\phi_n\}$, where*

$$\phi_{n+1} \triangleq \phi[v_n], \qquad 0 \le n < \infty \quad \text{and} \quad \phi_\infty \triangleq \phi[v_\infty] \tag{4.21}$$

*are the unique strictly positive roots of the functions $(Gv_n)(\cdot)$, $0 \le n \le \infty$, as in (4.17). Moreover,*

$$v_{n+1}(\cdot) = (Hv_n)(\cdot), \qquad 0 \le n < \infty, \quad \text{and} \quad v_\infty(\cdot) = (Hv_\infty)(\cdot). \tag{4.22}$$

*For every $n \ge 0$, we have $v_n(\phi_n) = v_n'(\phi_n) = 0$, and $v_{n+1}(\cdot)$ and $v_n(\cdot)$ satisfy*

$$\left\{ \begin{array}{r} v_{n+1}(\phi) < 0 \\ [\mathcal{A}_0 - (\lambda+\lambda_0)]v_{n+1}(\phi) + (g + \lambda_0(Kv_n))(\phi) = 0 \end{array} \right\}, \qquad \phi \in (0, \phi_{n+1}), \tag{4.23}$$

$$\left\{ \begin{array}{r} v_{n+1}(\phi) = 0 \\ [\mathcal{A}_0 - (\lambda+\lambda_0)]v_{n+1}(\phi) + (g + \lambda_0(Kv_n))(\phi) > 0 \end{array} \right\}, \qquad \phi \in (\phi_{n+1}, \infty). \tag{4.24}$$

*The function $v_\infty(\cdot)$ satisfies $v_\infty(\phi_\infty) = v_\infty'(\phi_\infty) = 0$ and the variational inequalities*

$$\left\{ \begin{array}{r} v_\infty(\phi) < 0 \\ [\mathcal{A}_0 - (\lambda+\lambda_0)]v_\infty(\phi) + (g + \lambda_0(Kv_\infty))(\phi) = 0 \end{array} \right\}, \qquad \phi \in (0, \phi_\infty), \tag{4.25}$$

$$\left\{ \begin{array}{r} v_\infty(\phi) = 0 \\ [\mathcal{A}_0 - (\lambda+\lambda_0)]v_\infty(\phi) + (g + \lambda_0(Kv_\infty))(\phi) > 0 \end{array} \right\}, \qquad \phi \in (\phi_\infty, \infty). \tag{4.26}$$

PROOF. Since $v_0(\cdot) \equiv 0$ is continuous, $v_1(\cdot) \triangleq (Jv_0)(\cdot) = (Hv_0)(\cdot)$ by (3.2) and Proposition 4.1, and $v_1(\cdot)$ is continuous by Lemma 4.3. Then an induction on $n$ and repeated applications of (3.2), Proposition 4.1 and Lemma 4.3 prove that every $v_n(\cdot)$, $0 \le n < \infty$, is continuous, and that the equalities on the left in (4.22) hold. Since $v_\infty(\cdot)$ is the uniform pointwise limit of the sequence $\{v_n(\cdot)\}_{n\ge 0}$ of continuous functions on $\mathbb{R}_+$ by Lemma 3.2, it is also continuous. Therefore, Lemma 3.1 and Proposition 4.1 also imply that $v_\infty(\cdot) = (Jv_\infty)(\cdot) = (Hv_\infty)(\cdot)$, which is the second equality in (4.22). The remainder of the corollary now follows from (4.22) and Lemma 4.3. □

PROPOSITION 4.2. *The pointwise limit $v_\infty(\cdot)$ of the sequence $\{v_n(\cdot)\}_{n\ge 0}$ in (3.2) and the value function $V(\cdot)$ of the optimal stopping problem in (2.6) coincide. The first entrance time $\widetilde{\tau}_{[\phi_\infty,\infty)}$ of the process $\Phi$ of (1.3)–(2.9) into the half interval $[\phi_\infty,\infty)$ is optimal for the Bayesian sequential change detection problem in (1.2) and (2.5).*

PROOF. Let $\widetilde{\tau}$ be an $\mathbb{F}$-stopping time, and $\widetilde{\tau}_{l,r} \triangleq \widetilde{\tau}_{[0,l]} \wedge \widetilde{\tau}_{[r,\infty)}$ for some $0 < l < r < \infty$. Then for every $\phi > 0$, the chain rule implies that $e^{-\lambda(\widetilde{\tau}\wedge\widetilde{\tau}_{l,r})} v_\infty(\Phi_{\widetilde{\tau}\wedge\widetilde{\tau}_{l,r}})$ equals

$$\begin{aligned}
v_\infty(\Phi_0) &+ \int_0^{\widetilde{\tau}\wedge\widetilde{\tau}_{l,r}} e^{-\lambda t}([\mathcal{A}_0 - (\lambda+\lambda_0)]v_\infty(\Phi_{t-}) + \lambda_0 (Kv_\infty)(\Phi_{t-}))\,dt \\
&+ \int_0^{\widetilde{\tau}\wedge\widetilde{\tau}_{l,r}} e^{-\lambda t}\mu\Phi_{t-} v'_\infty(\Phi_{t-})\,dX_t \\
&+ \int_0^{\widetilde{\tau}\wedge\widetilde{\tau}_{l,r}} \int_E e^{-\lambda t}\left[v_\infty\left(\frac{\lambda_1}{\lambda_0} f(z)\Phi_{t-}\right) - v_\infty(\Phi_{t-})\right] q(dt,dz)
\end{aligned}$$

in terms of the $(\mathbb{P}_0^\phi, \mathbb{F})$-compensated counting measure $q(dt,dz) \triangleq p(dt,dz) - \lambda_0\,dt\,\nu_0(dz)$ on $[0,\infty)\times E$. The stochastic integrals with respect to $X$ and $q$ are square-integrable martingales stopped at some $\mathbb{F}$-stopping time with finite expectation by Remark 3.2, since continuous $v'_\infty(\cdot)$ is bounded on $[l,r]$, and $v_\infty(\cdot)$ is bounded everywhere. Therefore, taking expectations of both sides implies that $\mathbb{E}_0^\phi[e^{-\lambda(\widetilde{\tau}\wedge\widetilde{\tau}_{l,r})} v_\infty(\Phi_{\widetilde{\tau}\wedge\widetilde{\tau}_{l,r}})]$ equals

$$\begin{aligned}
v_\infty(\phi) &+ \mathbb{E}_0^\phi\left[\int_0^{\widetilde{\tau}\wedge\widetilde{\tau}_{l,r}} e^{-\lambda t}([\mathcal{A}_0 - (\lambda+\lambda_0)]v_\infty(\Phi_{t-}) + \lambda_0(Kv_\infty)(\Phi_{t-}))\,dt\right] \\
&\geq v_\infty(\phi) - \mathbb{E}_0^\phi\left[\int_0^{\widetilde{\tau}\wedge\widetilde{\tau}_{l,r}} e^{-\lambda t} g(\Phi_{t-})\,dt\right] \\
&= v_\infty(\phi) - \mathbb{E}_0^\phi\left[\int_0^{\widetilde{\tau}\wedge\widetilde{\tau}_{l,r}} e^{-\lambda t} g(\Phi_t)\,dt\right],
\end{aligned}$$

since $[\mathcal{A}_0 - (\lambda+\lambda_0)]v_\infty(\cdot) + (g + \lambda_0(Kv_\infty))(\cdot) \geq 0$ because of the variational inequalities in (4.25) and (4.26). Since $v_\infty(\cdot)$ is bounded and continuous, letting $l \to 0$, $r \to \infty$, the bounded and monotone convergence theorems give

$$\mathbb{E}_0^\phi[e^{-\lambda\widetilde{\tau}} v_\infty(\Phi_{\widetilde{\tau}})] \geq v_\infty(\phi) - \mathbb{E}_0^\phi\left[\int_0^{\widetilde{\tau}} e^{-\lambda t} g(\Phi_t)\,dt\right]$$

and

$$\mathbb{E}_0^\phi\left[\int_0^{\widetilde{\tau}} e^{-\lambda t} g(\Phi_t)\,dt\right] \geq v_\infty(\phi)$$

for every $\mathbb{F}$-stopping time $\widetilde{\tau}$, because $v_\infty(\cdot)$ is nonpositive. By taking the infimum of both sides of the second inequality over all $\widetilde{\tau} \in \mathbb{F}$, we find that $V(\phi) \geq v_\infty(\phi)$ for every $\phi \in (0,\infty)$.

If we replace every $\widetilde{\tau}$ above with the $\mathbb{P}_0^\phi$-a.s. finite (by Lemma 4.1) $\mathbb{F}$-stopping time $\widetilde{\tau}_{[\phi_\infty,\infty)}$, then we have $\mathbb{E}_0^\phi[e^{-\lambda\widetilde{\tau}_{[\phi_\infty,\infty)}} v_\infty(\Phi_{\widetilde{\tau}_{[\phi_\infty,\infty)}})] = 0$ and every inequality becomes an equality by Corollary 4.1. Therefore, $V(\phi) \le \mathbb{E}_0^\phi[\int_0^{\widetilde{\tau}_{[\phi_\infty,\infty)}} e^{-\lambda t} g(\Phi_t)\,dt] = v_\infty(\phi)$. Hence,

$$V(\phi) = \mathbb{E}_0^\phi\left[\int_0^{\widetilde{\tau}_{[\phi_\infty,\infty)}} e^{-\lambda t} g(\Phi_t)\,dt\right] = v_\infty(\phi) \qquad \text{for every } \phi > 0, \tag{4.27}$$

and $\tau_{[\phi_\infty,\infty)}$ is optimal for (2.6) for every $\phi > 0$.

The same equalities at $\phi = 0$ and optimality of the stopping time $\tau_{[\phi_\infty,\infty)}$ when the initial state is 0 follow after taking limits in (4.27) as $\phi$ goes to zero if we prove that three functions in (4.27) are continuous at $\phi = 0$. The function $v_\infty(\cdot)$ is continuous on $[0,\infty)$ by Corollary 4.1. If we let $\widetilde{\tau} = \widetilde{\tau}_{[\phi_\infty,\infty)}$ and $\tau = \tau_{[\phi_\infty,\infty)}$ as in (4.7), then the strong Markov property of $\Phi$ at the first jump time $T_1$ gives

$$\begin{aligned}
w(\phi) &\triangleq \mathbb{E}_0^\phi\left[\int_0^{\widetilde{\tau}_{[\phi_\infty,\infty)}} e^{-\lambda t} g(\Phi_t)\,dt\right] \\
&= \mathbb{E}_0^\phi\left[\int_0^{\widetilde{\tau}\wedge T_1} e^{-\lambda t} g(\Phi_t)\,dt + \mathbf{1}_{\{\widetilde{\tau} > T_1\}} \int_{T_1}^{\widetilde{\tau}} e^{-\lambda t} g(\Phi_t)\,dt\right] \\
&= \mathbb{E}_0^\phi\left[\int_0^{\tau} e^{-\lambda t} \mathbf{1}_{\{t<T_1\}} g(Y_t^{\Phi_0})\,dt + \mathbf{1}_{\{\tau > T_1\}} e^{-\lambda T_1} w\left(\frac{\lambda_1}{\lambda_0} f(Z_1) Y_{T_1}^{\Phi_0}\right)\right] \\
&= \mathbb{E}_0^\phi\left[\int_0^{\tau} e^{-(\lambda+\lambda_0)t} (g + \lambda_0 (Kw))(Y_t^{\Phi_0})\,dt\right] = (H_{\phi_\infty} w)(\phi), \qquad \phi \ge 0,
\end{aligned}$$

which is continuous at $\phi = 0$ by (4.14). It remains to show that $\phi \mapsto V(\phi)$ is continuous at $\phi = 0$. Let us denote by $\widetilde{\tau}_h$ and $\tau_h$ the stopping times $\widetilde{\tau}_{[h,\infty)}$ and $\tau_{[h,\infty)}$ for every $h > 0$, as in (4.7). Since $g(\phi) < 0$ for $0 \le \phi < \lambda/c$, it is never optimal to stop before $\Phi$ reaches $[\lambda/c,\infty)$, and for every $0 < h \le \lambda/c$, we have

$$\begin{aligned}
V(0) &= \mathbb{E}_0^0\left[\int_0^{\widetilde{\tau}_h\wedge T_1} e^{-\lambda t} g(\Phi_t)\,dt + e^{-\lambda(\widetilde{\tau}_h \wedge T_1)} V(\Phi_{\widetilde{\tau}_h\wedge T_1})\right] \\
&= \mathbb{E}_0^0\left[\int_0^{\tau_h\wedge T_1} e^{-\lambda t} g(\Phi_t)\,dt + e^{-\lambda(\tau_h \wedge T_1)} V(\Phi_{\tau_h\wedge T_1})\right] \\
&= \mathbb{E}_0^0\left[\int_0^{\tau_h} e^{-\lambda t} \mathbf{1}_{\{t<T_1\}} g(Y_t^{\Phi_0})\,dt + \mathbf{1}_{\{\tau_h<T_1\}} e^{-\lambda\tau_h} V(Y_{\tau_h}^{\Phi_0})\right. \\
&\qquad\qquad \left. + \mathbf{1}_{\{\tau_h \ge T_1\}} e^{-\lambda T_1} V\left(\frac{\lambda_1}{\lambda_0} f(Z_1) Y_{T_1}^{\Phi_0}\right)\right] \\
&= \mathbb{E}_0^0\left[\int_0^{\tau_h} e^{-(\lambda+\lambda_0)t} g(Y_t^{\Phi_0})\,dt + e^{-(\lambda+\lambda_0)\tau_h} V(Y_{\tau_h}^{\Phi_0})\right. \\
&\qquad\qquad \left. + \int_0^{\tau_h} e^{-(\lambda+\lambda_0)t} \lambda_0 (KV)(Y_t^{\Phi_0})\,dt\right]
\end{aligned}$$

$$= \mathbb{E}_0^0\left[\int_0^{\tau_h} e^{-(\lambda+\lambda_0)t}(g+\lambda_0(KV))(Y_t^{\Phi_0})\,dt\right] + V(h)\frac{\psi(0)}{\psi(h)},$$

because $\{\widetilde{\tau}_h < T_1\} = \{\tau_h < T_1\}$ and $\widetilde{\tau}_h \wedge T_1 = \tau_h \wedge T_1$. Since $V(\cdot)$ is bounded, the first term on the right-hand side vanishes as $h \searrow 0$ by Remark 3.2. Because $\psi(0) > 0$, $\lim_{h\searrow 0} \psi(0)/\psi(h) = 1$. Therefore, $\lim_{h\searrow 0} V(h)$ exists and equals $V(0)$. Hence, $V(\phi)$ is also continuous at $\phi = 0$. $\square$

REMARK 4.1. The value function $V(\cdot) \equiv v_\infty(\cdot)$ can be approximated uniformly by the elements of the sequence $\{v_n(\cdot)\}$ at any desired level of accuracy according to the inequalities in (3.10). Since $\{v_n(\cdot)\}_{n\geq 0}$ decreases to $v_\infty(\cdot)$, the optimal continuation region

$$C \triangleq \{\phi \geq 0 : V(\phi) < 0\} \equiv \{\phi \geq 0 : v_\infty(\phi) < 0\} = [0, \phi_\infty)$$

is the increasing limit of $C_n \triangleq \{\phi \geq 0 : v_n(\phi) < 0\} = [0, \phi_n)$, $n \geq 0$, and $\phi_\infty = \lim_{n\to\infty} \uparrow \phi_n$.

Moreover, for every $\varepsilon > 0$ and for every $n \geq 1$ such that $[\lambda_0/(\lambda+\lambda_0)]^n < c\varepsilon$, the stopping time $\widetilde{\tau}_{[\phi_n,\infty)} = \inf\{t \geq 0; \Phi_t \geq \phi_n\}$ is $\varepsilon$-optimal for (2.6). More precisely,

$$\begin{aligned} V(\phi) &\leq \mathbb{E}_0^\phi\left[\int_0^{\widetilde{\tau}_{[\phi_n,\infty)}} e^{-\lambda t} g(\Phi_t)\,dt\right] \\ &\leq V(\phi) + \frac{1}{c}\left(\frac{\lambda_0}{\lambda+\lambda_0}\right)^n, \qquad \phi \geq 0, n \geq 1. \end{aligned}$$

PROOF. We shall prove the last displayed equation. Since $\widetilde{\tau}_{[\phi_n,\infty)}$ is the $\mathbb{P}_0^\phi$-a.s. finite $\mathbb{F}$-stopping time by Lemma 4.1, as shown for $\widetilde{\tau}_{[\phi_\infty,\infty)}$ in the proof of Proposition 4.2, Itô's rule and the localization argument imply that

$$\begin{aligned} &\mathbb{E}_0^\phi[e^{-\lambda\widetilde{\tau}_{[\phi_n,\infty)}} v_\infty(\Phi_{\widetilde{\tau}_{[\phi_n,\infty)}})] - v_\infty(\phi) \\ &\quad = \mathbb{E}_0^\phi\left[\int_0^{\widetilde{\tau}_{[\phi_n,\infty)}} e^{-\lambda t}([\mathcal{A}_0 - (\lambda+\lambda_0)]v_\infty + \lambda_0(Kv_\infty))(\Phi_t)\,dt\right] \\ &\quad = -\mathbb{E}_0^\phi\left[\int_0^{\widetilde{\tau}_{[\phi_n,\infty)}} e^{-\lambda t} g(\Phi_t)\,dt\right], \end{aligned}$$

since $([\mathcal{A}_0 - (\lambda+\lambda_0)]v_\infty + g + \lambda_0(Kv_\infty))(\phi) = 0$ for every $\phi \in (0, \phi_n) \subseteq (0, \phi_\infty)$ according to (4.25). Therefore, for every $\phi \geq 0$, we have

$$\begin{aligned} v_\infty(\phi) &= \mathbb{E}_0^\phi\left[\int_0^{\widetilde{\tau}_{[\phi_n,\infty)}} e^{-\lambda t} g(\Phi_t)\,dt\right] + \mathbb{E}_0^\phi[e^{-\lambda\widetilde{\tau}_{[\phi_n,\infty)}} v_\infty(\Phi_{\widetilde{\tau}_{[\phi_n,\infty)}})] \\ &\geq \mathbb{E}_0^\phi\left[\int_0^{\widetilde{\tau}_{[\phi_n,\infty)}} e^{-\lambda t} g(\Phi_t)\,dt\right] + \mathbb{E}_0^\phi[e^{-\lambda\widetilde{\tau}_{[\phi_n,\infty)}} v_n(\Phi_{\widetilde{\tau}_{[\phi_n,\infty)}})] \\ &\quad - \frac{1}{c}\left(\frac{\lambda_0}{\lambda+\lambda_0}\right)^n \end{aligned}$$

by the second inequality in (3.10). The result now follows immediately because we have $\mathbb{P}_0^{\phi}$-a.s. $v_n(\Phi_{\widetilde{\tau}_{[\phi_n,\infty)}}) = 0$ by Lemma 4.1 and Corollary 4.1. $\square$

## 5. Quickest detection of a simultaneous change in several independent Wiener and compound Poisson processes.

Suppose that at the disorder time $\Theta$, the drift of a *d-dimensional Wiener process* $\vec{W} = (W^{(1)}, \ldots, W^{(d)})$ changes from $\vec{0}$ to $\vec{\mu} = (\mu^{(1)}, \ldots, \mu^{(d)})$ for some $1 \le d < \infty$ and $\vec{\mu} \in \mathbb{R}^d \setminus \{\vec{0}\}$. Then in the model of Section 2, the likelihood-ratio process of (2.3) and its dynamics in (2.4) become

$$L_t \triangleq \exp\left\{\vec{\mu}\vec{X}_t - \left(\frac{\|\vec{\mu}\|^2}{2} + \lambda_1 - \lambda_0\right)t\right\} \prod_{n:0<T_n\le t} \left(\frac{\lambda_1}{\lambda_0} f(Z_n)\right), \qquad t \ge 0,$$

$$dL_t = L_t \vec{\mu}\, d\vec{X}_t + L_{t-} \int_E \left(\frac{\lambda_1}{\lambda_0} f(z) - 1\right)[p(dt \times dz) - \lambda_0\, dt\, \nu_0(dz)], \qquad t \ge 0,$$

in terms of the $d$-dimensional observation process

$$\vec{X}_t = \vec{X}_0 + \vec{\mu}(t-\Theta)^+ + \vec{W}_t, \qquad t \ge 0. \tag{5.1}$$

The representation in (2.5) of the minimum Bayes risk $U(\cdot)$ in terms of the value function $V(\cdot)$ of the optimal stopping problem in (2.6) for the conditional odds-ratio process $\Phi$ of (1.3) remains valid, but instead of (2.8), the dynamics of $\Phi$ now become

$$\begin{aligned} d\Phi_t &= (\lambda + a\Phi_t)\, dt + \Phi_t \vec{\mu}\, d\vec{X}_t \\ &\quad + \Phi_{t-} \int_E \left(\frac{\lambda_1}{\lambda_0} f(z) - 1\right) p(dt \times dz), \qquad t \ge 0. \end{aligned} \tag{5.2}$$

However, if we define

$$\mu \triangleq \|\vec{\mu}\| = \sqrt{(\mu^{(1)})^2 + \cdots + (\mu^{(d)})^2} \quad \text{and} \quad X \triangleq \frac{1}{\|\vec{\mu}\|}\vec{\mu}\vec{X}, \tag{5.3}$$

then the one-dimensional process $X$ is a $(\mathbb{P}_0^{\phi}, \mathbb{F})$-Wiener process $\mathbb{P}_0^{\phi}$-independent of the marked point process $(T_n, Z_n)_{n\ge1}$. In terms of the Wiener process $X$ and the new scalar $\mu \ne 0$ in (5.3), the dynamics in (5.2) of the sufficient statistic $\Phi$ can be rewritten exactly as in (2.8). Hence, quickest detection of a change from $\vec{0}$ to $\vec{\mu}$ in the drift of a multi-dimensional Wiener process is equivalent to quickest detection of a change from 0 to $\mu \equiv \|\vec{\mu}\|$ in the scalar drift of a suitable one-dimensional Wiener process. This is true both in the absence and presence of an independent and observable marked point process whose local characteristics change at the same time $\Theta$ as described earlier.

Suppose that, in addition to the process $\vec{X}$ in (5.1), $m$ compound Poisson processes $(T_n^{(i)}, Z_n^{(i)})_{n\geq 1}$, $1 \leq i \leq m$, independent of each other and the process $\vec{W}$, are observed on some common mark space $(E, \mathcal{E})$. At the same disorder time $\Theta$, their arrival time and mark distribution change from $(\lambda_0^{(i)}, \nu_0^{(i)})$ to $(\lambda_1^{(i)}, \nu_1^{(i)})$, respectively. Then their superposition forms a new marked point process $(T_n, Z_n)_{n\geq 1}$, which is independent of $\vec{W}$, and whose local characteristics are

$$(5.4) \qquad (\lambda_0, \nu_0) \triangleq \left( \sum_{i=1}^m \lambda_0^{(i)}, \sum_{i=1}^m \frac{\lambda_0^{(i)}}{\lambda_0} \nu_0^{(i)} \right) \quad \text{and}$$

$$(\lambda_1, \nu_1) \triangleq \left( \sum_{i=1}^m \lambda_1^{(i)}, \sum_{i=1}^m \frac{\lambda_1^{(i)}}{\lambda_1} \nu_1^{(i)} \right)$$

before and after the disorder time, respectively. Therefore, the solution method of the previous section, as summarized by Remark 4.1, can be applied directly with the new choices in (5.3) and (5.4) of parameters $\mu, \lambda_0, \lambda_1$, probability distributions $\nu_0, \nu_1$ on $(E, \mathcal{E})$ and processes $X$ and $(T_n, Z_n)_{n\geq 1}$.

**6. Numerical examples.** We describe briefly the numerical computation of the fundamental solution $\psi(\cdot)$ of the ODE $(\lambda + \lambda_0)h = \mathcal{A}_0 h$ and successive approximations $v_n(\cdot)$, $n \geq 0$, in (3.2) of the value function $V(\cdot)$ in (2.6). These computations are based on Kushner and Dupuis's [13] Markov chain approximation and Monte Carlo estimation of certain expectations. We use these methods on several examples and illustrate that reduction in Bayes risk can be significant if multiple sources are used simultaneously in order to detect an observable disorder time.

6.1. *Calculation of the function $\psi(\cdot)$ over a fine grid on some interval $[0, z]$.* Let us fix a number $z > 0$, grid size $h > 0$ and an integer $N$ such that $z = Nh$. Denote by $S^h$ the collection of grid points $z_n = nh$, $n \geq 0$. Set $\psi(z) = 1$ (or to any other positive constant). Then we can calculate the function $\psi(\cdot)$ on the grid $S^h$ according to (4.9) if we can evaluate

$$(6.1) \qquad \mathbb{E}_0^y[\exp\{-(\lambda + \lambda_0)\tau_z\}] \qquad \text{for every } y, z > 0.$$

To do that, we will approximate the diffusion $Y$ in (2.11) with a continuous-time process $\{\xi^h(t); t \geq 0\}$ obtained from a discrete-time Markov chain $\{\xi_n^h; n \geq 0\}$ on the state space $S^h$ by replacing unit-length sojourn times with state-dependent deterministic times. The derivation of one-step transition probabilities $p^h(y, v)$, $y, v \in S^h$ of the Markov chain $\{\xi_n^h; n \geq 0\}$ and "interpolation intervals" $\Delta t^h(y)$, $y \in S^h$ become more transparent if we set our goal to approximate the more general expectation

$$(6.2) \qquad V_\beta(y) \triangleq \mathbb{E}_0^y \left[ \int_0^{\tau_z} e^{-\beta t} k(Y_t)\, dt + e^{-\beta \tau_z} g(Y_{\tau_z}) \right], \qquad 0 < y < z,$$

for some fixed $z \in S^h$, discount rate $\beta \geq 0$, and bounded functions $k(\cdot)$ and $g(\cdot)$. Let us study $V_0$ first (namely, $\beta = 0$). If we denote the drift and diffusion coefficients of the process $Y$ by $\mu(\cdot)$ and $\sigma(\cdot)$, then, under certain regularity conditions, we expect $V_0(\cdot)$ to solve the second order ODE

$$\frac{\sigma^2(y)}{2} V_0''(y) + \mu(y) V_0'(y) + k(y) = 0, \qquad 0 < y < z,$$

subject to boundary condition $V_0(z) = g(z)$. If we replace $V_0''(y)$ and $V_0'(y)$ with their finite-difference approximations

$$\frac{V_0^h(y+h) + V_0^h(y-h) - 2V_0^h(y)}{h^2},$$

$$\frac{V_0^h(y+h) - V_0^h(y)}{h} \mathbf{1}_{[0,\infty)}(\mu(y)) + \frac{V_0^h(y) - V_0^h(y-h)}{h} \mathbf{1}_{(-\infty,0)}(\mu(y)),$$

respectively, then we obtain

$$\begin{aligned} &\frac{\sigma^2(y)}{2} \frac{V_0^h(y+h) + V_0^h(y-h) - 2V_0^h(y)}{h^2} \\ &\qquad + \frac{V_0^h(y+h) - V_0^h(y)}{h} \mu^+(y) \\ &\qquad + \frac{V_0^h(y) - V_0^h(y-h)}{h} \mu^-(y) + k(y) = 0, \qquad 0 < y < z. \end{aligned}$$

Rearranging the terms implies that $V_0^h(y)$ equals

$$\begin{aligned} &V_0^h(y-h) \frac{(\sigma^2(y)/2) + h\mu^-(y)}{\sigma^2(y) + h|\mu(y)|} + V_0^h(y+h) \frac{(\sigma^2(y)/2) + h\mu^+(y)}{\sigma^2(y) + h|\mu(y)|} \\ &\qquad + \frac{h^2}{\sigma^2(y) + h|\mu(y)|} k(y), \end{aligned}$$

which can be rewritten as

$$\begin{aligned} V_0^h(y) &= V_0^h(y-h) p^h(y, y-h) \\ &\qquad + V_0^h(y+h) p^h(y, y+h) + \Delta t^h(y) k(y) = 0, \end{aligned} \tag{6.3}$$

for every $y \in S^h \cap [0, z]$, if we define

$$\left\{ \begin{aligned} p^h(y, y \pm h) &\triangleq \frac{(\sigma^2(y)/2) + h\mu^\pm(y)}{\sigma^2(y) + h|\mu(y)|} \\ \Delta t^h(y) &\triangleq \frac{h^2}{\sigma^2(y) + h|\mu(y)|} \end{aligned} \right\}, \qquad y \in S^h. \tag{6.4}$$

Let $\{\xi_n^h; n \geq 0\}$ be the discrete-time Markov chain on $S^h$ with transition probabilities $p^h(y, y \pm h)$, $y \in S^h$, in (6.4), and define the continuous-time process $\{\xi^h(t); t \geq 0\}$ on the same space by adding the "interpolation interval" $\Delta t^h(\xi_n^h)$ before the jump from $\xi_n^h$ to $\xi_{n+1}^h$, namely,

$$\xi^h(t) \triangleq \xi_n^h \qquad t \in [t_n^h, t_{n+1}^h), n \geq 0$$
$$\text{where } t_0^h \equiv 0,\ t_{n+1}^h \triangleq t_n^h + \Delta t_n^h, n \geq 0 \text{ and } \Delta t_n^h \triangleq \Delta t^h(\xi_n^h)$$

are deterministic functions of the embedded discrete-time Markov chain $(\xi_n^h)_{n\geq 0}$. Then the solution $V_0^h(y)$, $y \in S^h \cap [0,z]$ of (6.3) with the boundary condition $V_0^h(z) = 0$ is the same as the expectation

$$V_0^h(y) = \mathbb{E}_0^y\left[\int_0^{\tau^h} k(\xi^h(t))\,dt + g(\xi^h(\tau^h))\right], \qquad y \in S^h \cap [0,z]. \tag{6.5}$$

The process $\{\xi^h(t); t \geq 0\}$ is *locally consistent* with $\{Y_t; t \geq 0\}$; and therefore, that process and the function $V_0^h(\cdot)$ well approximate $\{Y_t; t \geq 0\}$ and $V_0(\cdot)$, respectively; see Kushner and Dupuis [13] for the details. In general,

$$V_\beta^h(y) \triangleq \mathbb{E}_0^y\left[\int_0^{\tau^h} e^{-\beta t} k(\xi^h(t))\,dt + e^{-\beta\tau^h} g(\xi^h(\tau^h))\right], \qquad y \in S^h \cap [0,z], \tag{6.6}$$

is a good approximation of the function $V_\beta(\cdot)$ in (6.2), and if we define

$$N^h \triangleq \inf\{n \geq 0 : \xi_n^h = z\},$$

then (6.6) simplifies to

$$V_\beta^h(y) = \mathbb{E}_0^y\left[\sum_{n=0}^{N^h-1} k(\xi_n^h) e^{-\beta t_n^h} \frac{1 - e^{-\beta \Delta t_n^h}}{\beta} + \exp\{-\beta t_{N^h}^h\} g(z)\right], \qquad y \in S^h \cap [0,z]. \tag{6.7}$$

In (6.1), $\beta = \lambda + \lambda_0$, $k \equiv 0$, and $g \equiv 1$. Thus, (6.1) is approximated well by

$$\mathbb{E}_0^y[\exp\{-(\lambda + \lambda_0) t_{N^h}^h\}] \qquad \text{for } y \in S^h \cap [0,z] \text{ as well as } y \in S^h \cap [z,\infty). \tag{6.8}$$

Finally, we can estimate (6.8) by using Monte Carlo simulation in the following way:

(i) Set the initial state $\xi_0^h = y$.

(ii) Simulate the Markov chain $\xi_n^h$ until the first time $N^h$ that it hits the state $z \in S^h$.

(iii) Calculate $\exp\{-(\lambda + \lambda_0) \sum_{n=0}^{N^h-1} \Delta t^h(\xi_n^h)\}$, which is now a sample estimate of (6.8).

(iv) Repeat until the standard error of the sample average of individual estimates obtained from independent simulation runs reduces to an acceptable level. Report upon stopping the sample average as the approximate value of (6.8).

For the calculations in (4.9), notice that initial state $y$ and target state $z$ are always adjacent. This usually helps to keep the number of simulations low. In the detection problem, the dynamics in (6.4) of the Markov chain that approximates the diffusion $Y$ in (2.11) become

$$\left\{ \begin{aligned} p^h(y, y \pm h) &= \frac{(\mu^2/2)y^2 + h(\lambda + ay)^{\pm}}{\mu^2 y^2 + h|\lambda + ay|} \\ \Delta t^h(y) &= \frac{h^2}{\mu^2 y^2 + h|\lambda + ay|} \end{aligned} \right\}, \qquad y \in S^h. \tag{6.9}$$

We choose $h$ so small that $p^h(h, 2h) \gg p^h(h, 0)$, that is, reaching to 0 from inside $S^h$ is made almost impossible.

6.2. *Calculation of the successive approximations $v_n(\cdot)$, $n \geq 0$, in (4.22) of the value function $V(\cdot)$ in (2.6).* Recall from (3.2), Corollary 4.1 and (4.22) that bounded, nonpositive and nondecreasing functions $v_n(\cdot)$, $n \geq 0$, can be found by successive applications of the operator $H$ in (4.13) and (4.18). Therefore, it is enough to describe the calculation of $(Hw)(\cdot)$ for a bounded, nonpositive and nondecreasing function $w(\cdot)$.

Since the function $\psi(\cdot)$ is now available, the unique root $\phi[w]$ of $(Gw)(\phi) = 0$ in (4.17) can be found by solving numerically the equation

$$\int_0^{\phi[w]} z^{2[(a/\mu^2)-1]} e^{-2\lambda/(\mu^2 z)} \psi(z)[g(z) + \lambda_0 (Kw)(z)]\, dz = 0.$$

By Lemma 4.3, we have $(Hw)(\phi) = 0$ for every $\phi \geq \phi[w]$. Let $S^h$ denote once again the grid points $z_n = nh$, $n < N$, where $h > 0$ is small and $z_N = \phi[w]$. Then by simulating the approximate Markov chain $\{\xi_n^h; n \geq 0\}$ with transition probabilities and interpolation interval given in (6.9), we can approximate $(Hw)(\phi)$ on $S^h$ with the Monte Carlo estimate of

$$\mathbb{E}_0^{z_n}\left[\sum_{n=0}^{N^h-1} (g + \lambda_0 (Kw))(\xi_n^h) e^{-(\lambda+\lambda_0)t_n^h} \frac{1 - e^{-(\lambda+\lambda_0)\Delta t_n^h}}{\lambda + \lambda_0}\right] \qquad \text{at every } z_n \in S^h; \tag{6.10}$$

compare (6.2) and (6.7) with (4.13) and (6.10) when $r = \phi[w]$.

**Initialization**: Calculate simultaneously

- the increasing fundamental solution $\psi(\cdot)$ by simulating (6.8) on the interval $\phi \in [0, u]$,
- the function

$$I(z) \triangleq z^{2[(a/\mu^2)-1]} e^{-2\lambda/(\mu^2 z)} \psi(z) \propto \frac{2\psi(z)}{\sigma^2(z)B(z)} \qquad \text{for every } z \in [0, u],$$

- and the unique strictly positive solution $u \triangleq \phi[-1/c]$ of the inequality

$$0 \le \int_0^u I(z)\Big(g(z) - \frac{\lambda_0}{c}\Big)\, dz \propto (Gw)(u), \qquad \text{where } w(\cdot) \equiv -\frac{1}{c}.$$

For example, set initially $u = 2\lambda/c$, calculate the functions $\psi(\phi)$ and $I(\phi)$ for every $\phi \in [0, u]$. If the above inequality is satisfied, then stop, otherwise, double $u$ and repeat. Since $-1/c \le v_n(\cdot) \le 0$ for every $n \ge 0$ by Remark 3.3, we have $0 \le \phi_n \le u$ for every $n \ge 0$ by Lemma 4.2. Recall that $g(z) = z - (\lambda/c)$, $z \ge 0$, and $(Gw)(\cdot)$ are as in (2.6) and (4.17), respectively. Finally, set $n = 0$, $\phi_n = \lambda/c$, and $v_n(\phi) = 0$ for every $0 \le \phi \le u$.

**Step 1**: Calculate the function $(Kv_n)(\cdot)$ by using (3.4) and unique root $\widehat{\phi}_n \equiv \phi_\ell[v_n]$ of the increasing function $(g + \lambda_0(Kv_n))(\cdot)$.

**Step 2**: Find the unique strictly positive solution $r = \phi_{n+1}$ of the equation

$$0 = \int_0^r I(z)(g + (Kv_n))(z)\, dz \propto (Gv_n)(r).$$

The solution $\phi_{n+1}$ is located in the interval $(\widehat{\phi}_n \vee \phi_n, u)$, and Newton's method may be used to find it.

**Step 3**: Set $v_{n+1}(\phi) = 0$ for every $\phi_{n+1} \le \phi \le u$, and find $v_{n+1}(\phi)$ for every $0 \le \phi \le \phi_{n+1}$ by simulating (6.10). Increase $n$ by one and go to **Step 1**.

FIG. 1. *An algorithm that calculates the approximation $v_n(\cdot)$ of the value function $V(\cdot)$ and the critical thresholds $\phi_n$ for every $n \ge 0$; see Remark 4.1.*

6.3. *Examples.* Figure 1 describes an algorithm that calculates the approximations $v_n(\cdot)$, $n \ge 0$, of the value function $V(\cdot)$ by means of the tools described in Sections 6.1 and 6.2. In the following examples, we employ that algorithm to compute the approximations $v_n(\cdot)$, $n \ge 0$, until the maximum difference between two successive functions is reduced to an acceptable level. The termination of the algorithm with guaranteed error bounds follows from Lemma 3.2, which also provides an upper bound on the number of successive approximations.

Nine panels in Figure 2 display the approximate value functions corresponding to nine examples. In each example, the observation process is $(X, N)$; the process $X$ is a one-dimensional Wiener process that gains a drift $\mu$ after the disorder time $\Theta$, and $N$ is a simple Poisson process whose arrival rate changes from $\lambda_0$ to $\lambda_1$ at time $\Theta$. In all of the nine examples, we have $c = 1$ and $\lambda = \lambda_1 = 1$ [see (1.2) and (2.1)]; however, the post-disorder drift $\mu$ of $X$ and the pre-disorder arrival rate $\lambda_0$ of $N$ are different. Across every row, $\mu$ increases while $\lambda_0$ does not change. Across every column, $\lambda_0$ increases while $\mu$ does not change.

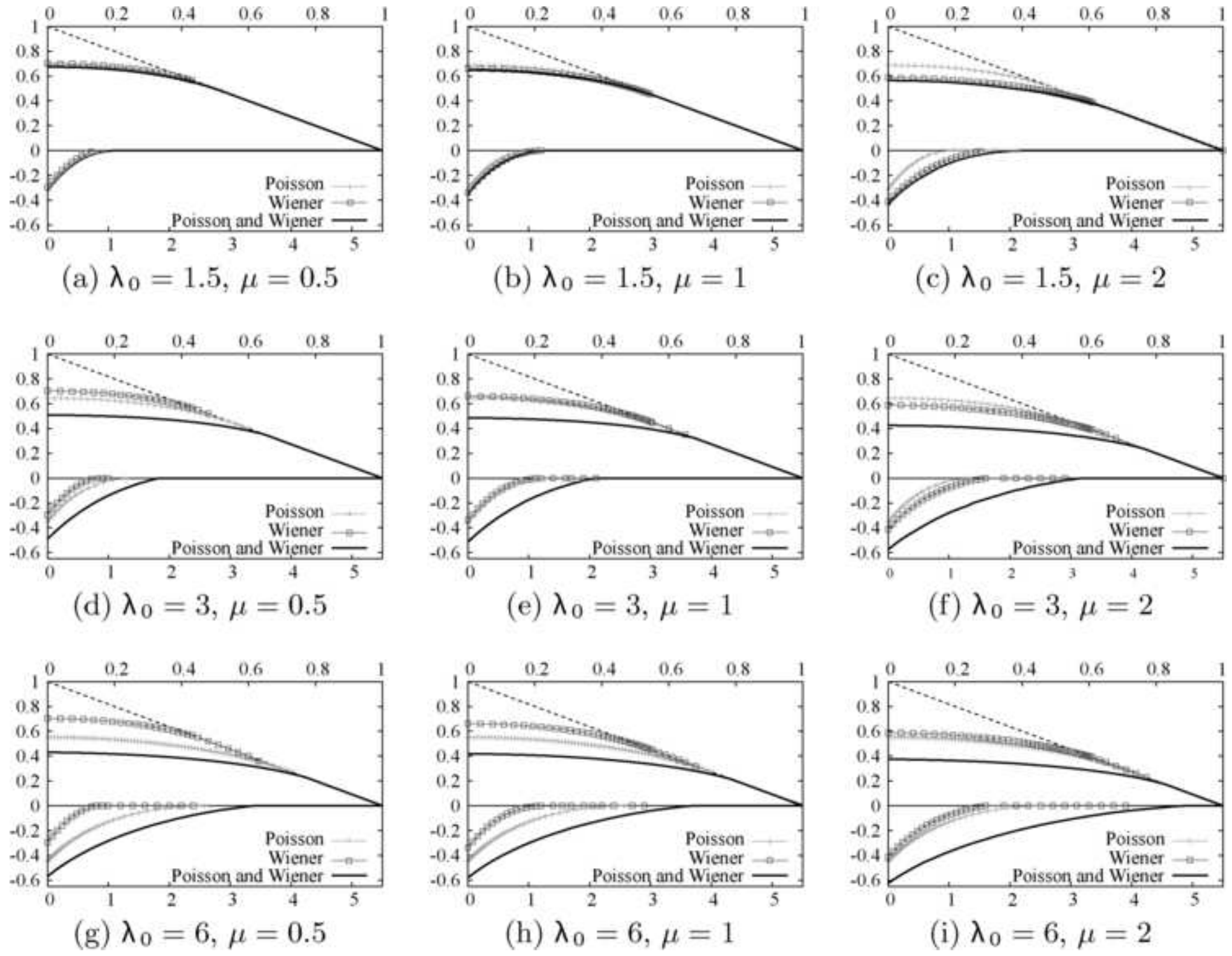

FIG. 2. *The solutions of the sequential disorder-detection problems for different pre-disorder arrival rate $\lambda_0$ of the Poisson process and post-disorder drift $\mu$ of the Wiener process. In each case, $\lambda_1 = \lambda = c = 1$. The upper part of each panel displays the Bayes risks $U(\cdot)$, $U_p(\cdot)$, $U_X(\cdot)$, and the lower part displays the value functions $V(\cdot)$, $V_p(\cdot)$, $V_X(\cdot)$ of the corresponding optimal stopping problems. Solid curves are the functions $U(\cdot)$ and $V(\cdot)$. Curves with "+" are the optimal Bayes risk $U_p(\cdot)$ and the value function $V_p(\cdot)$ if only the Poisson process is observable, and curves with "□" are the Bayes risk $U_X(\cdot)$ and the value function $V_X(\cdot)$ if only the Wiener process is observable. The dashed line in the upper part is the mapping $\pi \mapsto 1 - \pi$. It is optimal to raise a disorder alarm as soon as the process $\Phi/(1+\Phi)$ of (2.3) enters into the region where $U(\pi) = 1 - \pi$; equivalently, as soon as $\Phi$ enters the region where $V(\phi) = 0$.*

The graph in each panel is divided in two parts. The upper part shows the optimal Bayes risk $U(\cdot)$ of (2.5) on $[0, 1]$ displayed on the upper horizontal axis, and the lower part shows the value function $V(\cdot)$ of the stopping problem in (2.6) on $\mathbb{R}_+$ displayed on the lower horizontal axis. Both $U(\cdot)$ and $V(\cdot)$ are plotted with solid curves. We compare those functions with $U_p(\cdot)$, $V_p(\cdot)$, $U_X(\cdot)$ and $V_X(\cdot)$, where $U_p(\cdot)$ and $U_X(\cdot)$ are obtained by taking the infimum in (2.5) over the stopping times of (smaller) natural filtrations $\mathbb{F}^p$ and $\mathbb{F}^X$ of $N$ and $X$, respectively. On the other hand, $V_p(\cdot)$ and $V_X(\cdot)$ are the value functions of the optimal stopping problems analogous to (2.6), that is,

$$V_p(\phi) \triangleq \inf_{\tau \in \mathbb{F}^p} \mathbb{E}_0^\phi \left[ \int_0^\tau e^{-\lambda t} \left( \Phi_t^{(p)} - \frac{\lambda}{c} \right) dt \right],$$
$$V_X(\phi) \triangleq \inf_{\tau \in \mathbb{F}^X} \mathbb{E}_0^\phi \left[ \int_0^\tau e^{-\lambda t} \left( \Phi_t^{(X)} - \frac{\lambda}{c} \right) dt \right],$$

where

$$\Phi_t^{(p)} \triangleq \frac{\mathbb{P}\{\Theta \le t \mid \mathcal{F}_t^p\}}{\mathbb{P}\{\Theta > t \mid \mathcal{F}_t^p\}} \quad \text{and} \quad \Phi_t^{(X)} \triangleq \frac{\mathbb{P}\{\Theta \le t \mid \mathcal{F}_t^X\}}{\mathbb{P}\{\Theta > t \mid \mathcal{F}_t^X\}};$$

$U_p(\cdot)$, $V_p(\cdot)$ and $U_X(\cdot)$, $V_X(\cdot)$ are related to each other in the same way as $U(\cdot)$, $V(\cdot)$ are in (2.5).

The differences in the Bayes risks $U_p(\cdot)$, $U_X(\cdot)$ and $U(\cdot)$ provide insights about the contributions of observing the processes $X$ and $N$ separately or simultaneously to the efforts of detecting the disorder time $\Theta$. Sometimes, the Poisson process provides more information than the Wiener process, as in (d), (g) and (h); sometimes, the Wiener process provides more information than the Poisson, as in (b), (c) and (f); and some other times, the difference is negligible, as in (a), (e) and (i). In every case, observing the Poisson *and* Wiener processes at the same time provides more information, which is often *significantly larger* than two processes can provide separately, as in (i), (e), (f), (h), (d) and (g).

Intuitively, we expect the contributions to increase as $\mu$ and $\lambda_0$ are pulled farther apart from 0 and $\lambda_1$, respectively. The examples displayed in Figure 2 are consistent with this expectation. The Bayes risks $U_N(\cdot)$ and $U(\cdot)$ are shifting downward across every column, and $U_X(\cdot)$ and $U(\cdot)$ do the same across every row.

In (a), $\mu$ and $\lambda_0$ are relatively close to 0 and $\lambda_1$, respectively; therefore, observing both processes at the same time does not improve the optimal Bayes risk. Observing only one of them will thus reduce costs without increasing risks. As the post-disorder drift $\mu$ of $X$ is increased along the first row, both $U_X(\cdot)$ and $U(\cdot)$ improve gradually. The function $U_X(\cdot)$ stays close to $U(\cdot)$ because the process $X$ provides more information than $N$ for the detection of the disorder time. Especially in (c), one may choose not to observe the process $N$ anymore in order to lower the observation costs. Similarly, if $\mu$ is close to 0, an increase in the difference between $\lambda_0$ and $\lambda_1$ makes $U_p(\cdot)$ drive $U(\cdot)$ to lower levels; see the first column.

6.4. *Numerical comparison with Baron and Tartakovsky's asymptotic analysis.* Let us denote the Bayes risk $R_\tau(\pi)$ in (1.2), minimum Bayes risk $U(\pi)$ in (2.5) by $R_\tau(\pi, c)$ and $U(\pi, c)$, respectively, in order to display explicitly their dependence on the cost $c$ per unit detection delay. Let us also define

$$\begin{aligned} \phi(c) &\triangleq \frac{(\mu^2/2) + \lambda_0 + \lambda_1[\log(\lambda_1/\lambda_0) - 1] + \lambda}{c} \quad \text{and} \\ f(c) &\triangleq -\frac{\log c}{\phi(c)}, \qquad c > 0. \end{aligned} \tag{6.11}$$

Baron and Tartakovsky ([1], Theorem 3.5) have shown that the stopping time $\tau(c) \triangleq \inf\{t \geq 0; \Phi_t \geq \phi(c)\}$ is asymptotically optimal and that the minimum Bayes risk $U(\pi, c)$ asymptotically equals $f(c)$ for every fixed $\pi \in [0, 1)$, as the detection delay cost $c$ decreases to zero, in the sense that

$$\lim_{c \searrow 0} \frac{U(\pi, c)}{f(c)} = \lim_{c \searrow 0} \frac{R_{\tau(c)}(\pi, c)}{f(c)} = 1 \qquad \text{for every } \pi \in [0, 1).$$

In this subsection we revisit the example displayed in Figure 2(h), where $\lambda_0 = 6$ and $\lambda_1 = \lambda = \mu = 1$. We have calculated optimal thresholds, minimum Bayes risks and their asymptotic expansions in (6.11) for eighteen values of $c$ $(0.02, 0.04, \ldots, 0.18, 0.20, 0.30, \ldots, 0.90, 1)$; see Figure 3. If only the Poisson or Wiener process is observable, then the asymptotic expansions of the optimal thresholds and their minimum Bayes risks $U_p(\cdot)$, $U_X(\cdot)$ also follow from (6.11) by setting $\mu = 0$ in the Poisson case and by letting $\lambda_0 = \lambda_1$ in the Wiener case, respectively. The critical thresholds and minimum Bayes risks are calculated in Figure 3(c) by using the numerical algorithm in Figure 1, in Figure 3(a) by using Dayanik and Sezer's [9], Figure 2 numerical algorithm, and in Figure 3(b) by solving numerically the integral equation

$$\int_0^{\phi_\infty} \frac{[w - (1/c)]\psi_X(w)}{e^{2/w}}\, dw = 0$$

for the critical value $\phi_\infty$ and by numerically evaluating $V_X(\phi) =$

$$\psi_X(\phi) \int_\phi^{\phi_\infty} \frac{2[w - (1/c)]\eta_X(w)}{e^{2/w}}\, dw + \eta_X(\phi) \int_0^\phi \frac{2[w - (1/c)]\psi_X(w)}{e^{2/w}}\, dw$$

in $U_X(\pi) = 1 - \pi + c(1 - \pi)V_X(\pi/[1 - \pi])$, in terms of

$$\psi_X(\phi) = 1 + \phi \quad \text{and} \quad \eta_X(\phi) = (1 + \phi) \int_\phi^\infty \frac{e^{2/w}}{w^2(1 + w)^2}\, dw;$$

see also Shiryaev ([17], page 201, Theorem 9).

Optimal critical thresholds and their asymptotic expansions seem to be in good agreement; this is especially clear for small $c$ values as Baron and Tartakovsky [1] predicted (as $c$ decreases, the distance between any two curves in the first row does not grow faster than the critical thresholds themselves, hence, the relative error converges to zero). In the second row, the Bayes risks at three fixed values, $\pi = 0, 0.5, 0.8$ (one in the middle and two close to end-points of the range $[0, 1]$), also seem in good agreement with the asymptotic expansions for small values of detection delay cost $c$. As a reference, we have also plotted the minimum Bayes risks at optimal critical thresholds, which do not have to agree closely with the asymptotic expansions, because in this case minimum Bayes risks are evaluated at different $\pi$ values as $c$ changes, and their asymptotics do not immediately fall inside the scope of Theorem 3.5 of Baron and Tartakovsky [1].

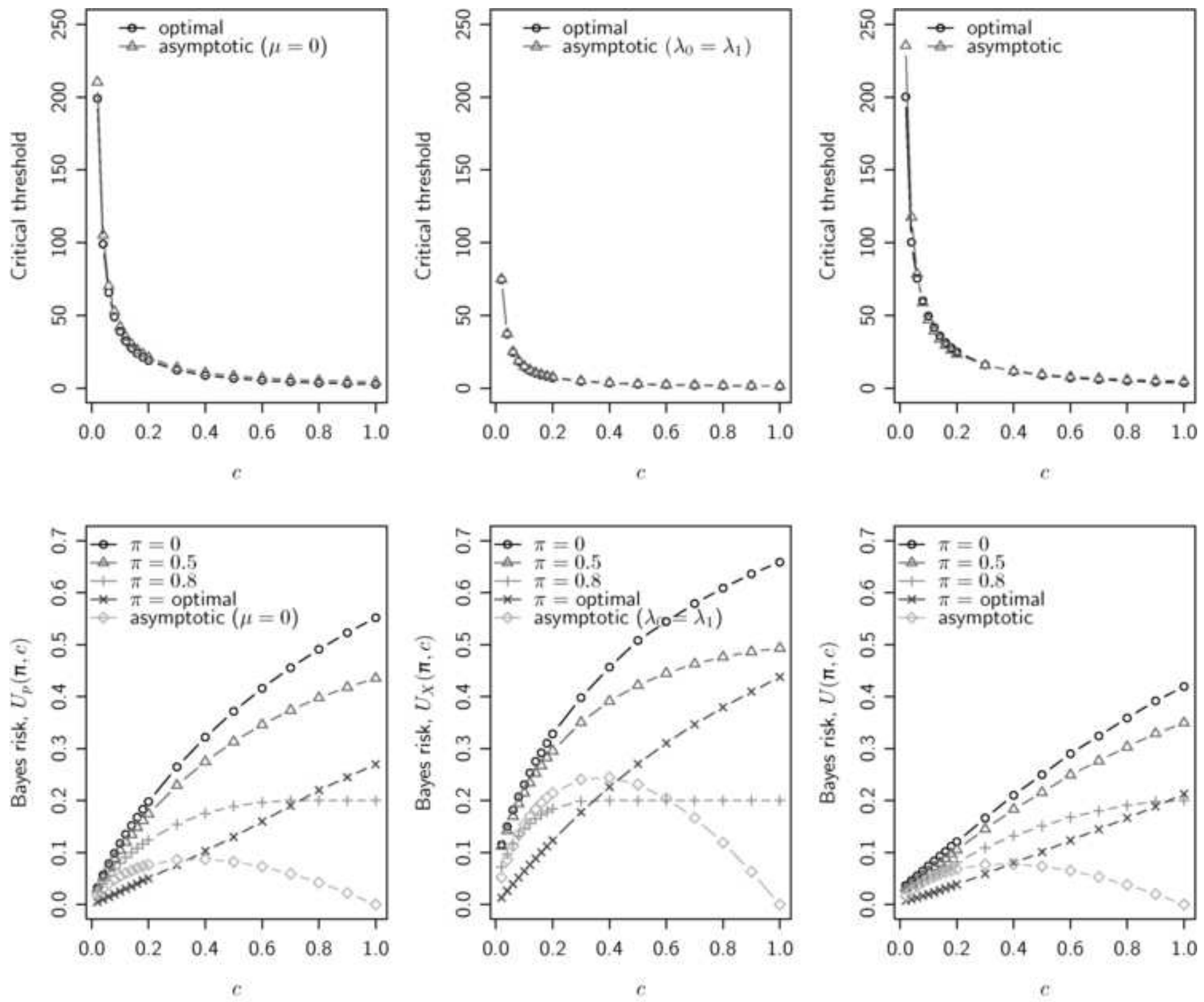

(a) Poisson (b) Wiener (c) Poisson and Wiener

FIG. 3. *Optimal thresholds, minimum Bayes risks and their asymptotic expansions at* $c = 0.02, 0.04, \ldots, 0.18, 0.20, 0.30, \ldots, 0.90, 1$ *for* (a) *Poisson,* (b) *Wiener and* (c) *combination of Poisson and Wiener observations* $(\lambda_0 = 6, \lambda_1 = \lambda = \mu = 1)$.

## APPENDIX

**A.1. The boundary behavior of the diffusion process $Y^y$.** Once we verify (4.2) and (4.3), the conclusions follow from Karlin and Taylor ([12],

Chapter 15), who expressed the quantities in (4.2) and (4.3) in terms of the measures $S(0,x] = \int_{0+}^{x} S(dy)$ and $M(0,x] = \int_{0+}^{x} M(dy)$, and integrals $\Sigma(0) = \int_{0+}^{x} S(0,\xi]M(d\xi)$, $N(0) = \int_{0+}^{x} M(0,\xi]S(d\xi)$ for the left boundary at 0, and $\Sigma(\infty) = \int_{x}^{\infty} S(\xi,\infty)M(d\xi)$, $N(\infty) = \int_{x}^{\infty} M(\xi,\infty)S(d\xi)$ for the right boundary at $\infty$. Since only the finiteness of $\Sigma(\cdot)$ and $N(\cdot)$ matters, the value of $x > 0$ in the domain of those integrals can be arbitrary. One finds that

$$S(dy) = c_1 y^{-2a/\mu^2} e^{2\lambda/(\mu^2 y)}\, dy$$

and

$$M(dy) = c_2 y^{2[(a/\mu^2)-1]} e^{-2\lambda/(\mu^2 y)}\, dy, \qquad y > 0;$$

above, as well as below, $c_1, c_2, \ldots$ will denote positive proportionality constants. Therefore, changing the integrating variable by setting $z = 1/y$ gives

$$\begin{aligned} S(x) - S(0+) &= \int_{0+}^{x} S(dy) \\ &= c_1 \int_{1/x}^{\infty} z^{(2a/\mu^2)-2} e^{(2\lambda/\mu^2)z}\, dz = +\infty \qquad \forall x > 0, \end{aligned}$$

and the first equality in (4.2) follows. After applying the same change of variable twice, the double integral in the same equation becomes

$$N(0) = c_3 \int_{1/x}^{\infty} \left( \int_{v}^{\infty} u^{\alpha} e^{-\beta u}\, du \right) v^{-\alpha-2} e^{\beta v}\, dv \tag{A.1}$$

in terms of $\alpha \triangleq -2a/\mu^2 \in \mathbb{R}$ and $\beta \triangleq 2\lambda/\mu^2 > 0$. Integrating the inner integral by parts $k \ge 0$ times gives that, for every $k \ge 0$,

$$\begin{aligned} \int_{v}^{\infty} u^{\alpha} e^{-\beta u}\, du = {} & \sum_{j=0}^{k-1} \frac{\alpha! \beta^{-(j+1)}}{(\alpha - j)!} v^{\alpha - j} e^{-\beta v} \\ & + \frac{\alpha! \beta^{-(k+1)}}{(\alpha - k)!} \int_{v}^{\infty} \beta u^{\alpha - k} e^{-\beta u}\, du. \end{aligned}$$

If $k \ge \alpha$, then $u \mapsto u^{\alpha-k}$ is decreasing and the integral on the right is less than or equal to $v^{\alpha-k} \int_{v}^{\infty} \beta e^{-\beta u}\, du = v^{\alpha-k} e^{-\beta v}$. Therefore,

$$\int_{v}^{\infty} u^{\alpha} e^{-\beta u}\, du \le \sum_{j=0}^{k} \frac{\alpha! \beta^{-(j+1)}}{(\alpha - j)!} v^{\alpha - j} e^{-\beta v}, \qquad k \ge \max\{0, \alpha\}.$$

Using this estimate in (A.1) implies that, for every $x > 0$,

$$\begin{aligned} N(0) &\le \int_{1/x}^{\infty} \left( \sum_{j=0}^{k} \frac{\alpha! \beta^{-(j+1)}}{(\alpha - j)!} v^{\alpha - j} e^{-\beta v} \right) v^{-\alpha-2} e^{\beta v}\, dv \\ &= \sum_{j=0}^{k} \frac{\alpha! \beta^{-(j+1)}}{(\alpha - j)!} \int_{1/x}^{\infty} v^{-j-2}\, dv < \infty, \end{aligned}$$

which completes the proof of (4.2). Since $S(0+) = -\infty$ and $N(0) < \infty$, the left boundary at 0 is an entrance-not-exit boundary.

For the proof of (4.3), notice that change of variable by $u = 1/y$ gives for every $z > 0$ that

$$\int_z^\infty S(dy) = \int_0^{1/z} u^{-\alpha-2} e^{\beta u}\, du \geq \int_0^{1/z} u^{-\alpha-2}\, du$$
$$= \begin{cases} -(\alpha+1)^{-1} z^{\alpha+1}, & \alpha+1<0, \\ \infty, & \alpha+1 \geq 0. \end{cases}$$

If $\alpha + 1 \geq 0$, then clearly $\Sigma(\infty) = \int_x^\infty \int_z^\infty S(dy)M(dz) = \infty$ for every $x > 0$. If $\alpha + 1 < 0$, then for every $x > 0$ we also have

$$\Sigma(\infty) = \int_x^\infty \int_z^\infty S(dy)M(dz) \geq \int_x^\infty -(\alpha+1)^{-1} z^{\alpha+1} M(dz)$$
$$= -(\alpha+1)^{-1} c_2 \int_x^\infty z^{\alpha+1} z^{-\alpha-2} e^{-\beta/z}\, dz = c_4 \int_x^\infty z^{-1} e^{-\beta/z}\, dz$$
$$\geq c_4 e^{-\beta/x} \int_x^\infty z^{-1}\, dz = \infty,$$

and the first equality in (4.3) is proved. Similarly, changing variable by $v = 1/y$ gives

$$\int_z^\infty M(dy) = \int_0^{1/z} v^\alpha e^{-\beta v}\, dv \geq e^{-\beta/z} \int_0^{1/z} v^\alpha\, dv$$
$$= \begin{cases} (\alpha+1)^{-1} z^{\alpha+1} e^{-\beta/z}, & \alpha+1>0, \\ \infty, & \alpha+1 \leq 0. \end{cases}$$

If $\alpha + 1 \leq 0$, then clearly $N(\infty) = \int_x^\infty \int_z^\infty M(dy)S(dz) = \infty$ for every $x > 0$. If $\alpha + 1 > 0$, then for every $x > 0$ we also have

$$N(\infty) = \int_x^\infty \int_z^\infty M(dy)S(dz) \geq \int_x^\infty (\alpha+1)^{-1} z^{\alpha+1} e^{-\beta/z} S(dz)$$
$$= c_1 \int_x^\infty (\alpha+1)^{-1} z^{\alpha+1} e^{-\beta/z} z^\alpha e^{\beta/z}\, dz = c_5 \int_x^\infty z^{2\alpha+1}\, dz$$
$$= c_6 z^{2(\alpha+1)}\big|_{z=x}^{z=\infty} = \infty,$$

which completes the proof of (4.3). Because $\Sigma(\infty) = N(\infty) = \infty$, the right boundary at $\infty$ is a natural boundary.

**A.2. Continuity of $\phi \mapsto (H_r w)(\phi)$ at $\phi = 0$.** We shall prove the second equality in (4.14), namely, $(H_r w)(0) = \lim_{\phi \searrow 0} \lim_{l \searrow 0} (H_{l,r} w)(\phi) \equiv \lim_{\phi \searrow 0} (H_r w)(\phi)$, which implies along with the first equality in (4.14) that

$\phi \mapsto (H_r w)(\phi)$ is continuous at $\phi = 0$. For every $0 < h < r$,

$$\begin{aligned}
(H_r w)(0) &= \mathbb{E}_0^0\left[\int_0^{\tau_h} e^{-(\lambda+\lambda_0)t}(g+\lambda_0(Kw))(Y_t^{\Phi_0})\,dt\right.\\
&\qquad \left. + \int_{\tau_h}^{\tau_r} e^{-(\lambda+\lambda_0)t}(g+\lambda_0(Kw))(Y_t^{\Phi_0})\,dt\right]\\
&= \mathbb{E}_0^0\left[\int_0^{\tau_h} e^{-(\lambda+\lambda_0)t}(g+\lambda_0(Kw))(Y_t^{\Phi_0})\,dt\right.\\
&\qquad \left. + e^{-(\lambda+\lambda_0)\tau_h}(H_r w)(Y_{\tau_h}^{\Phi_0})\right]\\
&= \mathbb{E}_0^0\left[\int_0^{\tau_h} e^{-(\lambda+\lambda_0)t}(g+\lambda_0(Kw))(Y_t^{\Phi_0})\,dt\right]\\
&\quad + (H_r w)(h)\mathbb{E}_0^0 e^{-(\lambda+\lambda_0)\tau_h}\\
&= \mathbb{E}_0^0\left[\int_0^{\tau_h} e^{-(\lambda+\lambda_0)t}(g+\lambda_0(Kw))(Y_t^{\Phi_0})\,dt\right] + (H_r w)(h)\frac{\psi(0)}{\psi(h)}\\
&\xrightarrow{h\searrow 0} 0 + \lim_{h\searrow 0}(H_r w)(h)\cdot 1,
\end{aligned}$$

where the second equality follows from the strong Markov property of $Y^{\Phi_0}$ applied at the $\mathbb{F}$-stopping time $\tau_h = \inf\{t \ge 0; Y_t^{\Phi_0} = h\}$, and the fourth equality from (4.8). As $h \searrow 0$, $\mathbb{P}_0^0$-a.s. $\tau_h \searrow 0$ since 0 is an entrance-not-exit boundary, and the integral and its expectation in the last equation vanish by the bounded convergence theorem. Moreover, since $\psi(0) \equiv \psi(0+) > 0$ by (4.4), we have $\lim_{h\searrow 0}\psi(0)/\psi(h) = 1$. Therefore, $\lim_{h\searrow 0}(H_r w)(h)$ must exist, and taking limits of both sides in the last displayed equation completes the proof.

**A.3. Calculation of $(H_{l,r}w)(\cdot)$ in (4.15).** Let us denote the function on the right-hand side of (4.15) by $\widehat{H}w(\phi)$, $l \le \phi \le r$. It can be rewritten in the more familiar form

$$\widehat{H}w(\phi) = \int_l^r G_{l,r}(\phi, z)(g + \lambda_0(Kw))(z)\,dz, \qquad l \le \phi \le r,$$

by means of the Green function

$$G_{l,r}(\phi, z) = \frac{\psi_l(\phi \wedge z)\eta_r(\phi \vee z)}{\sigma^2(z)W_{l,r}(z)}, \qquad l \le \phi, z \le r,$$

for the second order ODE

$$\text{(A.2)} \qquad [\mathcal{A}_0 - (\lambda + \lambda_0)]H(\phi) = -(g + \lambda_0(Kw))(\phi),$$
$$l < \phi < r, \text{ with boundary conditions } H(l+) = H(r-) = 0.$$

Therefore, the continuous function $\widehat{H}w(\phi)$, $l \le \phi \le r$, is twice continuously differentiable on $(l,r)$ and solves the boundary value problem in (A.2). If $\tau_{l,r} \triangleq \tau_{[0,l]} \wedge \tau_{[r,\infty)}$, Itô's rule gives

$$
\begin{aligned}
e^{-(\lambda+\lambda_0)\tau_{l,r}} &\widehat{H}w(Y^{\Phi_0}_{\tau_{l,r}}) - \widehat{H}w(\Phi_0) \\
&= \int_0^{\tau_{l,r}} e^{-(\lambda+\lambda_0)t}[\mathcal{A}_0 - (\lambda+\lambda_0)]\widehat{H}w(Y^{\Phi_0}_t)\,dt \\
&\quad + \int_0^{\tau_{l,r}} e^{-(\lambda+\lambda_0)t}\sigma(Y^{\Phi_0}_t)\widehat{H}w'(Y^{\Phi_0}_t)\,dX_t \\
&= -\int_0^{\tau_{l,r}} e^{-(\lambda+\lambda_0)t}(g+\lambda_0(Kw))(Y^{\Phi_0}_t)\,dt \\
&\quad + \int_0^{\tau_{l,r}} e^{-(\lambda+\lambda_0)t}\sigma(Y^{\Phi_0}_t)\widehat{H}w'(Y^{\Phi_0}_t)\,dX_t,
\end{aligned}
$$

where $\mathbb{P}^{\phi}_0$ a.s. $\widehat{H}w(Y^{\Phi_0}_{\tau_{l,r}}) = 0$, since $\widehat{H}w(l) = \widehat{H}w(r) = 0$ and the first exit time $\tau_{l,r}$ of the regular diffusion $Y^{\Phi}$ from the closed bounded interval $[l,r] \subsetneqq [0,\infty)$ is always $\mathbb{P}^{\phi}_0$ a.s. finite. Moreover, the stochastic integral with respect to the $(\mathbb{P}_0,\mathbb{F})$-Wiener process $X$ on the right-hand side has zero expectation because the derivative $\widehat{H}w'(\phi)$, given by

$$
\begin{aligned}
\psi_l'(\phi) &\int_\phi^r \frac{2\eta_r(z)}{\sigma^2(z)W_{l,r}(z)}(g+\lambda_0(Kw))(z)\,dz \\
&+ \eta_r'(\phi)\int_l^\phi \frac{2\psi_l(z)}{\sigma^2(z)W_{l,r}(z)}(g+\lambda_0(Kw))(z)\,dz,
\end{aligned}
$$

of $\widehat{H}w(\phi)$, is bounded on $\phi \in [l,r]$. Therefore, taking expectations of both sides gives

$$
\widehat{H}w(\phi) = \mathbb{E}^{\phi}_0\left[\int_0^{\tau_{l,r}} e^{-(\lambda+\lambda_0)t}(g+\lambda_0(Kw))(Y^{\Phi_0}_t)\,dt\right] \equiv (H_r w)(\phi), \qquad l \le \phi \le r.
$$

**Acknowledgment.** The authors thank an anonymous referee for the remarks and suggestions which improved the presentation of this paper.

## REFERENCES

[1] Baron, M. and Tartakovsky, A. G. (2006). Asymptotic optimality of change-point detection schemes in general continuous-time models. *Sequential Anal.* **25** 257–296. MR2244607

[2] Bayraktar, E., Dayanik, S. and Karatzas, I. (2005). The standard Poisson disorder problem revisited. *Stochastic Process. Appl.* **115** 1437–1450. MR2158013

[3] Bayraktar, E., Dayanik, S. and Karatzas, I. (2006). Adaptive Poisson disorder problem. *Ann. Appl. Probab.* **16** 1190–1261. MR2260062

[4] Borodin, A. N. and Salminen, P. (2002). *Handbook of Brownian Motion—Facts and Formulae*, 2nd ed. Birkhäuser, Basel. MR1912205

[5] Brémaud, P. (1981). *Point Processes and Queues*. Springer, New York. MR0636252

[6] Byington, C. S. and Garga, A. K. (2001). *Handbook of Multisensor Data Fusion*. CRC Press, Boca Raton, FL.

[7] Cartea, I. and Figueroa, M. (2005). Pricing in electricity markets: A mean reverting jump diffusion model with seasonality. *Appl. Math. Finance* **12** 313–335.

[8] Çinlar, E. (2006). Jump-diffusions. Blackwell-Tapia Conference, 3–4 November 2006. Available at http://www.ima.umn.edu/2006-2007/SW11.3-4.06/abstracts.html#Cinlar-Erhan.

[9] Dayanik, S. and Sezer, S. O. (2006). Compound Poisson disorder problem. *Math. Oper. Res.* **31** 649–672. MR2281222

[10] Gapeev, P. V. (2005). The disorder problem for compound Poisson processes with exponential jumps. *Ann. Appl. Probab.* **15** 487–499. MR2115049

[11] Itô, K. and McKean, Jr., H. P. (1974). *Diffusion Processes and Their Sample Paths*. Springer, Berlin. MR0345224

[12] Karlin, S. and Taylor, H. M. (1981). *A Second Course in Stochastic Processes*. Academic Press, New York. MR0611513

[13] Kushner, H. J. and Dupuis, P. (2001). *Numerical Methods for Stochastic Control Problems in Continuous Time*, 2nd ed. Springer, New York. MR1800098

[14] Peskir, G. and Shiryaev, A. (2006). *Optimal Stopping and Free Boundary Problems*. Birkhäuser, Basel. MR2256030

[15] Peskir, G. and Shiryaev, A. N. (2002). Solving the Poisson disorder problem. In *Advances in Finance and Stochastics* (K. Sandmann and P. Schonbucher, eds.) 295–312. Springer, Berlin. MR1929384

[16] Polyanin, A. D. and Zaitsev, V. F. (2003). *Handbook of Exact Solutions for Ordinary Differential Equations*, 2nd ed. Chapman and Hall/CRC, Boca Raton, FL. MR2001201

[17] Shiryaev, A. N. (1978). *Optimal Stopping Rules*. Springer, New York. MR0468067

[18] Weron, R., Bierbrauer, M. and Trück, S. (2004). Modeling electricity prices: Jump diffusion and regime switching. *Phys. A Statist. Mech. Appl.* **336** 39–48.

S. Dayanik
Department of Operations Research
and Financial Engineering
and
The Bendheim Center for Finance
Princeton University
Princeton, New Jersey 08544
USA
E-mail: sdayanik@princeton.edu

H. V. Poor
School of Engineering
and Applied Science
Princeton University
Princeton, New Jersey 08544
USA
E-mail: poor@princeton.edu

S. O. Sezer
Department of Mathematics
University of Michigan
Ann Arbor, Michigan 48109
USA
E-mail: sezer@umich.edu